\documentclass[11pt, twoside]{article}

\usepackage[english]{babel}

\usepackage{amssymb}
\usepackage{mathrsfs}
\usepackage{amsmath}
\usepackage{amsthm}
\usepackage{amsfonts}
\usepackage{latexsym}
\usepackage{indentfirst}
\usepackage{color}
\usepackage{graphicx}
\usepackage{txfonts}
\usepackage{enumerate}

\usepackage[colorlinks=true,
linkcolor=blue,
citecolor=red,
urlcolor=magenta,
backref=page
]{hyperref}

\usepackage{txfonts}
\usepackage{anysize}

\allowdisplaybreaks

\newtheorem{theorem}{Theorem}[section]
\newtheorem{lemma}[theorem]{Lemma}

\newtheorem{proposition}[theorem]{Proposition}
\theoremstyle{definition}
\newtheorem{definition}[theorem]{Definition}
\newtheorem{remark}[theorem]{Remark}
\numberwithin{equation}{section}
\begin{document}
\title{\Large\bf Solvability of $L^2$-Neumann Problems for Fractional Laplace Equations on Bounded Lipschitz Domains
\footnotetext{\hspace{-0.35cm}2020 {\it Mathematics Subject Classification}. Primary 35R11; Secondary 35B65, 46E35.
\endgraf{\it Key words and phrases}. Fractional Laplacian, nonlocal Neumann condition, Lipschitz domain,
weighted gradient estimate, non-tangential maximal function.
\endgraf This project is partially supported by the National Natural Science Foundation of China (Grant No. 12431006)
and Longyuan Young Talents of Gansu Province.}}
\author{Wenxian Ma\,\footnote{Corresponding author, E-mail
\texttt{mawenxian@lzu.edu.cn}}\ \ and Sibei Yang}
\date{}
\maketitle
\vspace{-0.8cm}
\begin{center}
\begin{minipage}{13.5cm}\small
{\bf Abstract.} Let $n\geq2$, $s\in(0,1)$, and $\Omega\subset\mathbb{R}^n$
be a bounded Lipschitz domain. In this paper, we study the solvability of
the fractional Neumann problem
\begin{equation*}
\left\{\begin{aligned}
(-\Delta)^su&=0 & & \text{in}\ \ \Omega,\\
\mathcal{N}_s u&=g & & \text{in}\ \ \mathbb{R}^n\setminus\Omega,
\end{aligned}\right.
\end{equation*}
with $\mathcal{N}_s$ denoting the nonlocal normal derivative in $\mathbb{R}^n\setminus\Omega$ and $g:=a_sq$, where
$a_s(z):=(1-s)\frac{\delta(z)^{-s}}{1+\delta(z)^{n+s}}$ with
$\delta(z):=\mathrm{dist\,}(z,\partial\Omega)$ for any $z\in\mathbb{R}^n\setminus\Omega$,
and $q\in L^2(\mathbb{R}^n\setminus\Omega,a_s\,dz)$ satisfies
$\int_{\mathbb{R}^n\setminus\Omega}a_sq\,dz=0$. For the weak solution $u$
of this fractional Neumann problem normalized by
$\int_\Omega u\,dx=0$, we obtain the non-tangential maximal function estimate
\begin{equation*}
\left\|\left(\delta^{1-s}|\nabla u|\right)^*\right\|_{L^2(\mathbb{R}^n\setminus\Omega,\,a_s\,dz)}
\leq C\|q\|_{L^2(\mathbb{R}^n\setminus\Omega,\,a_s\,dz)},
\end{equation*}
where, for any measurable function $f$ on $\Omega$, $f^\ast$ denotes the non-tangential maximal
function of $f$ on $\mathbb{R}^n\setminus\Omega$. This result extends the classical $L^2$ solvability theory
for the Neumann problem of Laplace equations on bounded Lipschitz domains to the setting of
fractional Neumann problems, and also provides a Neumann counterpart to the
$L^2$-solvability result for fractional Dirichlet problems established
by Colombo, Fern\'andez-Real and Ros-Oton \cite{CFR26}.
In addition, we prove the fractional Sobolev estimate
$\|u\|_{H^{3s/2}(\Omega)}\leq C\|q\|_{L^2(\mathbb{R}^n\setminus\Omega,\,a_s\,dz)}$
and the weighted gradient estimate
\begin{equation*}
\left(\int_\Omega\delta^{2-3s}|\nabla u|^2\,dx\right)^{\frac12}
\leq C\|q\|_{L^2(\mathbb{R}^n\setminus\Omega,\,a_s\,dz)}.
\end{equation*}
To the best of our knowledge, for the inhomogeneous exterior nonlocal Neumann problem, the above $H^{3s/2}$
estimate, weighted gradient estimate, and non-tangential maximal estimate are even new on bounded smooth domains.
\end{minipage}
\end{center}
\vspace{-0.1cm}

\tableofcontents

\section{Introduction and main results}

Let $n\geq2$, $s\in(0,1)$, and $\Omega\subset\mathbb{R}^n$ be a bounded
Lipschitz domain. In this paper, we consider the fractional
Neumann problem
\begin{equation}\label{e1.1}
\left\{\begin{aligned}
(-\Delta)^su&=0 & & \text{in}\ \ \Omega,\\
\mathcal{N}_s u&=g & & \text{in}\ \ \mathbb{R}^n\setminus\Omega.
\end{aligned}\right.
\end{equation}
Recall that, for any
$u\in C_{\rm c}^\infty(\mathbb{R}^n)$ (the space of all infinitely
differentiable functions on $\mathbb R^n$ with compact support) and $x\in\mathbb{R}^n$,
\begin{equation*}
(-\Delta)^su(x):=c_{(n,s)}\,\mathrm{P.V.}\int_{\mathbb{R}^n}\frac{u(x)-u(y)}{|x-y|^{n+2s}}\,dy
\end{equation*}
and, for any $z\in\Omega^c:=\mathbb{R}^n\setminus\Omega$,
\begin{equation*}
\mathcal N_su(z):=c_{(n,s)}\int_\Omega\frac{u(z)-u(x)}{|z-x|^{n+2s}}\,dx
\end{equation*}
(see, for instance, \cite{DRV17,dglz12}), where
\begin{equation*}
c_{(n,s)}:=s2^{2s}\pi^{-\frac n2}\frac{\Gamma(\frac n2+s)}{\Gamma(1-s)},
\end{equation*}
$\Gamma$ denotes the Gamma function, and $\mathrm{P.V.}$
denotes the Cauchy principal value of the integral.
Motivated by the recent work of Colombo, Fern\'andez-Real,
and Ros-Oton \cite{CFR26} on the $L^2$ solvability of the fractional
Dirichlet problem on bounded Lipschitz domains, we study the $L^2$
solvability of the fractional Neumann problem \eqref{e1.1}.

For $x\in\mathbb{R}^n$, let $\delta(x):=\mathrm{dist\,}(x,\partial\Omega)$ with $\partial\Omega$
denoting the boundary of $\Omega$. Define the
\emph{exterior density} and \emph{measure} by setting, for any $z\in\Omega^c$,
\begin{equation*}
 a_s(z):=(1-s)\frac{\delta(z)^{-s}}{1+\delta(z)^{n+s}}
\qquad \text{and}\qquad
 d\sigma_s(z):=a_s(z)\,dz.
\end{equation*}
The exterior density and measure were introduced by Colombo, Fern\'andez-Real and Ros-Oton \cite{CFR26}
for the fractional Dirichlet problem, and an equivalent choice $\tau_s(z)=(1-s)\delta(z)^{-s}(1+\delta(z))^{-n-s}$
of the density was introduced by Grube and Hensiek \cite{GH24}. The \emph{exterior data space} of $g$ used
in this paper is
\begin{equation}\label{e1.3}
\mathcal X_s(\Omega^c):=L^2(\Omega^c,a_s^{-1}\,dz)
\end{equation}
with the norm
\begin{equation*}
\|g\|_{\mathcal X_s(\Omega^c)}:=\left[\int_{\Omega^c}\frac{|g(z)|^2}{a_s(z)}\,dz\right]^{1/2}.
\end{equation*}
Equivalently, $g=a_sq$ with $q\in L^2(\Omega^c,\sigma_s)$,
and their norms satisfy
\begin{equation}\label{e1.4}
\|g\|_{\mathcal X_s(\Omega^c)}=\|q\|_{L^2(\Omega^c,\,\sigma_s)}.
\end{equation}
In particular, $\mathcal X_s(\Omega^c)\subset L^1(\Omega^c)$,
since $\sigma_s(\Omega^c)<\infty$ by
Lemma~\ref{l2}\textup{(i)} below.

Let
\begin{equation}\label{e1.5}
\mathcal{Q}_\Omega:=\mathbb{R}^{2n}\setminus(\Omega^c\times\Omega^c)
\end{equation}
and
\begin{equation}\label{e1.6}
\mathcal{V}:=\left\{u:\mathbb{R}^n\to\mathbb{R}\text{ measurable}:
\|u\|_{L^2(\Omega)}^2+\mathcal{E}_\Omega(u,u)<\infty\right\},
\end{equation}
where
\begin{equation}\label{e1.7}
\mathcal{E}_\Omega(u,v):=\frac{c_{(n,s)}}2\iint_{\mathcal{Q}_\Omega}
\frac{[u(x)-u(y)][v(x)-v(y)]}{|x-y|^{n+2s}}\,dx\,dy.
\end{equation}
Given a function $g\in L^1(\Omega^c)$, define
\begin{equation}\label{e1.8}
\mathcal{V}_g:=\left\{v\in\mathcal{V}:\int_{\Omega^c}|v|^2|g|\,dz<\infty\right\}.
\end{equation}
We equip spaces $\mathcal{V}$ and $\mathcal{V}_g$ with the norms
\begin{equation*}
\|u\|_{\mathcal{V}}:=\left[\|u\|_{L^2(\Omega)}^2+\mathcal{E}_\Omega(u,u)\right]^{1/2}
\quad\text{and}\quad
\|v\|_{\mathcal{V}_g}:=\left[\|v\|_{\mathcal{V}}^2+\int_{\Omega^c}|v|^2|g|\,dz\right]^{1/2}
\end{equation*}
respectively. If $g\equiv0$, then $\mathcal{V}_g=\mathcal{V}$.

We now state the definition of the weak solution to problem \eqref{e1.1}.
\begin{definition}\label{def:weak}
Let $n\geq2$, $s\in(0,1)$, and $\Omega\subset\mathbb{R}^n$ be a bounded
Lipschitz domain. Assume that $g\in\mathcal X_s(\Omega^c)$ satisfies the
compatibility condition $\int_{\Omega^c}g\,dz=0$.
Then a function $u\in\mathcal{V}$ is called a \emph{weak solution} to problem \eqref{e1.1} if
$\int_\Omega u\,dx=0$ and, for any $v\in\mathcal{V}_g$,
\begin{equation}\label{e1.11}
\mathcal{E}_\Omega(u,v)=\int_{\Omega^c}gv\,dz.
\end{equation}
\end{definition}

\begin{remark}
We point out that, for the data space
\eqref{e1.3}, the identity \eqref{e1.11} also holds for
any test function $v\in\mathcal{V}$ (see Lemma~\ref{l3}).
For the present data space, it is not necessary to assume that
$g=\mathcal N_s\psi$ for some $C^2$ function $\psi$, as required
in \cite[Theorem 3.9]{DRV17} and \cite[Remark 1.4]{MY26}.
\end{remark}

\subsection{Main results}

Following \cite[Definition 4.5]{CFR26}, for any function $f$ on $\Omega$ and $z\in\Omega^c$, define
\begin{equation*}
\Gamma(z):=\{x\in\Omega:|x-z|\leq2\delta(x)\}
\quad\text{and}\quad
 f^*(z):=\sup_{x\in\Gamma(z)}|f(x)|,
\end{equation*}
with $f^*(z):=0$ if the cone $\Gamma(z)$ is empty. For $t\in(0,2)$ and an open set $\Omega\subset\mathbb{R}^n$,
denote by $H^t(\Omega)$ the restriction
type Sobolev space on $\Omega$ (see Section \ref{s2.1} for its detailed definition).

Now we state the main results of this paper.

\begin{theorem}\label{t1}
Let $n\geq2$, $s\in(0,1)$, and $\Omega\subset\mathbb{R}^n$ be a bounded
Lipschitz domain. Assume that $g\in\mathcal X_s(\Omega^c)$ satisfies the compatibility condition
$\int_{\Omega^c}g\,dz=0$. Then problem \eqref{e1.1} has a unique
weak solution $u\in\mathcal{V}$ in the sense of
Definition~\ref{def:weak}. Moreover, the solution $u$ belongs to $H^{3s/2}(\Omega)$
and satisfies
\begin{equation}\label{e1.12x}
\|u\|_{H^{3s/2}(\Omega)}\leq C\|g\|_{\mathcal X_s(\Omega^c)},
\end{equation}
\begin{equation}\label{e1.13}
\left[\int_\Omega\delta^{2-3s}|\nabla u|^2\,dx\right]^{1/2}
\leq C\|g\|_{\mathcal X_s(\Omega^c)},
\end{equation}
and
\begin{equation}\label{e1.14}
\left\|\left(\delta^{1-s}|\nabla u|\right)^*\right\|_{L^2(\Omega^c,\,\sigma_s)}
\leq C\|g\|_{\mathcal X_s(\Omega^c)},
\end{equation}
where the positive constant $C$ depends only on $n$, $s$, and $\Omega$.
\end{theorem}

\begin{remark}
In \cite{CFR26}, Colombo, Fern\'andez-Real and Ros-Oton established
the solvability of the fractional Dirichlet problem
\begin{equation*}
\left\{\begin{aligned}
(-\Delta)^su&=0 & & \text{in}\ \ \Omega,\\
u&=g & & \text{in}\ \ \Omega^c
\end{aligned}\right.
\end{equation*}
with $g\in L^2(\Omega^c,\,\sigma_s)$ and proved that
\begin{equation*}
\|u^*\|_{L^2(\Omega^c,\,\sigma_s)}
\leq C\|g\|_{L^2(\Omega^c,\,\sigma_s)}.
\end{equation*}
The non-tangential maximal function estimate \eqref{e1.14} is a Neumann counterpart of this result.
\end{remark}

\begin{remark}
Formally, as $s\to1^-$, the non-tangential maximal function estimate
\eqref{e1.14} reduces to the classical $L^2$ Neumann estimate of Jerison
and Kenig \cite{JK81}, with $\delta^{1-s}\to1$, the exterior measure
$\sigma_s$ converging weakly to surface measure on $\partial\Omega$
(see, for instance, \cite[Remark 1.4]{CFR26}), and the normalized exterior datum $q=g/a_s$
corresponding to the boundary Neumann datum.
\end{remark}

For a smooth function $u$ on $\Omega$ and an integer $m\geq1$, write $D^mu:=(\partial^\alpha u)_{|\alpha|=m}$.
Here $\alpha:=(\alpha_1,\ldots,\alpha_n)$ is a multi-index,
$|\alpha|:=\alpha_1+\cdots+\alpha_n$, and
$\partial^\alpha u:=\partial_1^{\alpha_1}\cdots\partial_n^{\alpha_n}u$.

The next result estimates these higher derivatives.

\begin{theorem}\label{t2}
Let $n\geq2$, $s\in(0,1)$, and $\Omega\subset\mathbb{R}^n$ be a bounded
Lipschitz domain. Assume that $g\in\mathcal X_s(\Omega^c)$ satisfies
the compatibility condition
$\int_{\Omega^c}g\,dz=0$ and $u\in\mathcal{V}$ is the
weak solution to problem \eqref{e1.1} in the sense of Definition~\ref{def:weak}.
Then, for any positive integer $m$,
\begin{equation*}
\left[\int_\Omega\delta^{2m-3s}|D^mu|^2\,dx\right]^{1/2}\leq C\|g\|_{\mathcal X_s(\Omega^c)},
\end{equation*}
where the positive constant $C$ depends only on $n$, $s$, $m$, and $\Omega$.
\end{theorem}

\subsection{Background and ideas of the proof}

Historically, the solvability of Dirichlet and Neumann problems of
elliptic partial differential equations has been
a central question in potential theory and partial differential
equations (see, for instance, \cite{Dah77,DK87,fjr78,JK81,JK81a,k94,KP93,KP95,Shen06,Shen07, Ve84}).
Passing from smooth domains to Lipschitz domains brings these problems into a setting that permits corners
and nonsmooth boundary graphs. In particular, in 1977, Dahlberg \cite{Dah77}
established reverse H\"older estimates for harmonic measure on Lipschitz
domains, which yields $L^2$ solvability of the Dirichlet problem for Laplace
equation on Lipschitz domains.
Subsequently, in 1981, Jerison and Kenig \cite{JK81} established the corresponding
$L^2$ Neumann theory. Verchota \cite{Ve84}
proved the invertibility of classical layer potential operators in the
$L^2$ setting on Lipschitz domains in 1984. In \cite{DK87}, Dahlberg and Kenig developed the solvability
theory for Neumann problem of Laplace equation on Lipschitz domains with general $L^p$ data and Hardy date.
Later developments included the Neumann problem for elliptic equations
with nonsmooth coefficients studied by Kenig and Pipher \cite{KP93,KP95}.
For higher-order elliptic systems with real constant coefficients on Lipschitz domains, Shen \cite{Shen06}
characterized $L^p$ solvability of the Dirichlet problem by a reverse H\"older inequality of order $p$.
Shen \cite{Shen07} later developed a layer-potential approach to the $L^p$ invertibility of boundary
operators, which yields solvability of Neumann problems for second-order elliptic systems and the biharmonic equation.

In recent years, nonlocal equations have become an increasingly active area of research in the field
of partial differential equations, with equations involving the fractional Laplacian serving as a
prototypical example (see, for instance, \cite{AFR23,BN23,cs07,CFR26,DRV17,DRVW25,dglz12,FK24,fls16,
jsw20,kms15,MY26,n21,n23,rs14a,RW26,rs14b,s07}). In particular, in \cite{DRV17}, Dipierro, Ros-Oton, and Valdinoci introduced
the fractional Neumann problem \eqref{e1.1} and developed its variational construction and probabilistic interpretations.
For inhomogeneous exterior Neumann data, Bersetche and
Borthagaray \cite{BB22} studied well-posedness of problem \eqref{e1.1} and finite element
approximation. Foghem and Kassmann \cite{FK24} developed a
general Hilbert space framework for the nonlocal Neumann problem. Grube and Hensiek \cite{GH24}
established robust trace and extension theorems for the function spaces adapted to
the nonlocal Neumann problem and corresponding energy
estimates for the solution. In particular, their exterior weight is equivalent to $a_s$.
Thus, the weighted data space and energy estimates have precedents.
Theorem~\ref{t1} in this paper concerns additional Sobolev regularity
and weighted control of derivatives for the compatible data space
$\mathcal X_s(\Omega^c)$.

Moreover, Audrito, Felipe-Navarro and Ros-Oton \cite{AFR23}
proved boundary H\"older regularity for weak solutions with homogeneous
exterior Neumann data. Dipierro, Ros-Oton, Valdinoci and Weidner
\cite{DRVW25} subsequently established optimal boundary regularity by
means of the Mellin transform, also with homogeneous exterior Neumann
data. More recently, Colombo, Fern\'andez-Real and Ros-Oton \cite{CFR26}
established a fractional analogue of Dahlberg's theorem and obtained
$L^2$ solvability of the exterior Dirichlet problem with nonlocal
non-tangential maximal function estimates. Extending this line of
research to the fractional Neumann problem requires control of
derivatives of the solution from boundary data prescribed throughout the exterior. In this case, both the
choice of exterior measure and the singular behavior near the boundary need to be considered.
The resulting estimates quantify how weighted $L^2$ exterior date controls
the Sobolev regularity of the solution and the non-tangential behavior of
its weighted gradient.

In \cite{MY26}, the authors of this paper established existence, uniqueness, pointwise
estimates, and Sobolev estimates for Neumann
functions of fractional Laplace equations on bounded Lipschitz domains.
The Sobolev estimates in
\cite[Theorem 1.7]{MY26} concern the Neumann function itself, while the
pointwise and gradient bounds in \cite[Theorems 1.5 and 1.8]{MY26}
involve an arbitrary positive loss $\eta$ in the singularity exponent.
Together with the representation formula, these bounds do not directly
yield the estimates in Theorem~\ref{t1} with the stated exponents and
only $\|g\|_{\mathcal X_s(\Omega^c)}$ on the right-hand side.
We also point out that the present problem is also distinct from
the local Neumann boundary conditions studied by Ros-Oton
and Weidner \cite{RW26}. The latter prescribe a normal derivative of
$u/\delta^{s-1}$ on $\partial\Omega$, whereas \eqref{e1.1} prescribes
the nonlocal normal derivative $\mathcal N_su$ throughout $\Omega^c$.

To obtain the main results of this paper, we first fix the values of $u$ in $\Omega$ and minimize the energy
$\mathcal{E}_\Omega(u,u)$ over all choices of exterior values. At the minimizer of $\mathcal E_\Omega(u,u)$,
the exterior value $Bu(z)$ is given, for any $z\in\Omega^c$, by
$$
Bu(z):=\frac{\int_\Omega |z-x|^{-n-2s}u(x)\,dx}{\int_\Omega |z-x|^{-n-2s}\,dx}.
$$
When the function $g$ in
\eqref{e1.1} is nonzero, problem \eqref{e1.1} determines the exterior value $u(z)$ as \begin{equation}\label{e1.15}
u(z)=Bu(z)+g(z)/(c_{(n,s)}\rho(z)),
\end{equation}
where $\rho(z)$ is the
integral of the kernel over $\Omega$, as defined in \eqref{e2.9}.
By this, in Lemma~\ref{l4}\textup{(iii)} we transform the original weak formulation \eqref{e1.11} to
\eqref{e2.26}. This transformation allows us to handle the equation and boundary conditions simultaneously in the domain.
The right-hand side of \eqref{e2.26} can be written as
$$
F_q(w)=\int_{\Omega^c} q\,Bw\,d\sigma_s.
$$
By H\"older's inequality and Hardy's inequality,
we obtain $F_q$ is bounded from $H^{s/2}(\Omega)$ to its dual $(H^{s/2}(\Omega))^*$.
Then, following the ideas of Savar\'e \cite{Sa98} and using the construction in \cite[Section 2.3]{BLN23},
we construct a smooth approximation to $u$ that takes into account the cone condition of the Lipschitz domain.
This, combined with an iteration argument and interpolation theory, we establish the Besov
and Sobolev regularities of the solution in Propositions \ref{t3} and \ref{t4}.

We then establish the weighted gradient estimate. When $3s/2<1$, the $H^{3s/2}(\Omega)$ regularity proved
above does not imply that the gradient is square integrable. We therefore use the Poisson representation of the
solution in the interior to estimate its derivative at a point in terms of differences of function values in
Section \ref{s5.1}. Since the equation is nonlocal, we must also handle the part away from that point; we
call this part the tail (see \eqref{e5.3}).
To control the part of the tail outside $\Omega$, we use the reflection map $R_j$ (see Section \ref{s4.3}) to
reflect an exterior point into \(\Omega\). We then split the difference between an exterior point $y$ and an
interior point $x$ into
$$
u(y)-u(x)=[u(y)-u(R_jy)]+[u(R_jy)-u(x)].
$$
The first term is handled using the characterization of exterior values in \eqref{e1.15}, the flattening map for the
Lipschitz domain, and the reflection difference estimate from Section \ref{s4.1}. The second term is treated similarly,
but does not require us to consider exterior values. Section \ref{s4} discusses the reflection interpolation technique
in detail. When $s<1/3$, the reflected difference can be estimated directly using Hardy's inequality. When $1/3<s<2/3$,
\cite[Proposition 3.8]{DNPV12} implies that the interior function has a boundary trace. We explicitly extend this
boundary value into the half-space. Subtracting this extension leaves the function with zero boundary value,
to which Hardy's inequality applies. When $s=1/3$, we combine estimates on either side of Sobolev order $1/2$ to
treat the same difference operator, including its weights. When $s=2/3$, the usual $H^1$ estimate suffices. When $s>2/3$,
we apply Hardy's inequality to the gradient, then use a change of coordinates to obtain the result on Lipschitz domains.
These arguments use the interpolation tools in \cite{Le23,BL76}, as well as Lemma~\ref{l1}. The reflected estimates
are proved in Lemmas~\ref{l8}, \ref{l9}, \ref{l10}, and \ref{l11}.

To estimate the non-tangential maximal function, we use the Whitney decomposition in Section~\ref{s5}. For a cube
$Q$ with side length $\ell_Q$, the set of exterior points whose associated cones intersect $Q$ has $\sigma_s$-measure
at most $C\ell_Q^{n-s}$. Thus, the contribution of $Q$ to the exterior non-tangential estimate is at most
$$
C\ell_Q^{n+2-3s}\sup_Q|\nabla u|^2.
$$
This term can be controlled using the interior gradient estimate in Section \ref{s5.2} and the tail estimate in
Section \ref{s5.1}. Summing over all the cubes then yields the desired conclusion.

The remainder of this paper is organized as follows.
In Section~\ref{s2}, we recall function spaces used in this paper and give
the desired energy estimates. In Section~\ref{s3},
we prove an energy comparison and the lifted Sobolev regularity estimate for the weak solution to \eqref{e1.1}.
In Section~\ref{s4}, we establish the exterior reflection difference bound for the weak solution to \eqref{e1.1}
by using an reflection argument. Finally, in Section~\ref{s5}, we establish an estimate for tails and reconstruct the
derivatives for the weak solution to \eqref{e1.1}, and further finish the proofs of Theorems \ref{t1} and \ref{t2}.

At the end of this section, we make some notational conventions. Let
${\mathbb N}:=\{1,2,\ldots\}$ and ${\mathbb N}_0:={\mathbb N}\cup\{0\}$.
We always denote by $C$ a positive constant which is independent of the main parameters involved, but it may
vary from line to line. We also use $C_{\alpha,\beta,\ldots}$ to denote a positive constant depending on
the indicated parameters $\alpha$, $\beta,\ldots$. The \emph{notation} $f\lesssim g$
means that $f\leq Cg$. If $f\lesssim g$ and $g\lesssim f$,we then write $f\sim g$. If $f\leq Cg$ and
$g=h$ or $g\leq h$, we then write $f\lesssim g=h$ or $f\lesssim g\leq h$. For any $x\in{\mathbb R}^n$ and $r\in(0,\infty)$,
we denote by $B(x,r):=\{y\in{\mathbb R}^n:\ |y-x|<r\}$ the ball with center $x$
and radius $r$. We use $\mathbf{1}_E$ to denote the
\emph{characteristic function} of a measurable set $E\subset{\mathbb R}^n$
and $\mathbf{0}$ to denote the \emph{origin} of ${\mathbb R}^n$.
For any $n\ge2$, let $\mathbb{R}^n_+:=\mathbb{R}^{n-1}\times(0,\infty)$, and
for any $x\in\mathbb{R}^n$, write $x:=(x',x_n)$ with $x'\in\mathbb{R}^{n-1}$.
The \emph{notation} $\Omega_0\Subset\Omega_1$ means that $\overline{\Omega}_0$ is a compact subset of
the open set $\Omega_1$, where $\overline{\Omega}_0$ denotes
the closure of $\Omega_0$ in $\mathbb{R}^n$.
For any real-valued functions $f,g\in L^2(\Omega)$, the inner product
$(f,g)_{L^2(\Omega)}$ is defined by $(f,g)_{L^2(\Omega)}:=\int_\Omega fg\,dx$.
For a Banach space $X$, denote by $X^*$ its dual space with the norm
\begin{equation*}
\|f\|_{X^*}:=\sup_{x\in X,\|x\|_X\leq1}|f(x)|.
\end{equation*}
For a bounded linear mapping $T:X\to Y$, denote by $T^*:Y^*\to X^*$
its adjoint mapping, defined by, for any $f\in Y^*$ and $x\in X$, $(T^*f)(x):=f(Tx)$.
We also use $I$ to denote the identity operator. Finally, let $K(x,y):=|x-y|^{-n-2s}$ for
any $x,\,y\in\mathbb{R}^n$ with $x\ne y$

\section{Preliminaries}\label{s2}
In this section, we recall Sobolev spaces, Besov spaces, and real interpolation spaces employed in
this paper and establish the energy estimates needed later.

\subsection{Function spaces}\label{s2.1}

\begin{definition}\label{d2.1}
Let $n\geq2$, and let $\Omega\subset\mathbb{R}^n$ be a bounded domain.
Then $\Omega$ is called a Lipschitz domain if there exist constants $r_0,L\in(0,\infty)$
such that, for any $x_0\in\partial\Omega$ and $r\in(0,r_0)$, there is a rigid change of coordinates for which
\begin{equation*}
\Omega\cap B(x_0,r)=\{x\in B(x_0,r):x_n>\varphi(x')\},
\end{equation*}
where $\varphi:\mathbb{R}^{n-1}\to\mathbb{R}$ is a Lipschitz function satisfying
$\|\nabla\varphi\|_{L^\infty(\mathbb{R}^{n-1})}\leq L$, and
$x':=(x_1,\ldots,x_{n-1})$.
\end{definition}

Let $n\ge1$. For any $f\in\mathcal S(\mathbb{R}^n)$, the \emph{Fourier transform}
$\mathcal Ff$ of $f$ is defined by setting, for any $\xi\in\mathbb{R}^n$,
\begin{equation*}
\mathcal Ff(\xi)=\widehat f(\xi):=
(2\pi)^{-n/2}\int_{\mathbb{R}^n}e^{-ix\cdot\xi}f(x)\,dx,
\end{equation*}
where $\mathcal S(\mathbb{R}^n)$ denotes the Schwartz space. The Fourier transform $\mathcal F$ can be extended by duality to
$\mathcal S'(\mathbb{R}^n)$, the space of tempered distributions, and
$\mathcal F^{-1}$ denotes its inverse. For any $t\in\mathbb R$,
the \emph{Sobolev space} $H^t(\mathbb{R}^n)$ is defined to be the set of all
$f\in\mathcal S'(\mathbb{R}^n)$ satisfying
$(1+|\cdot|^2)^{t/2}\widehat f\in L^2(\mathbb{R}^n)$, equipped with the norm
\begin{equation*}
\|f\|_{H^t(\mathbb{R}^n)}:=
\left[\int_{\mathbb{R}^n}(1+|\xi|^2)^t\left|\widehat f(\xi)\right|^2\,d\xi\right]^{1/2}.
\end{equation*}
In particular, $H^0(\mathbb{R}^n)=L^2(\mathbb{R}^n)$ with equivalent norms.

For an open set $\Omega\subset\mathbb{R}^n$, the restricted Sobolev space
$H^t(\Omega)$ is defined by restriction and is equipped with the norm
\begin{equation}\label{e2.1}
\|v\|_{H^t(\Omega)}:=\inf\left\{\|V\|_{H^t(\mathbb{R}^n)}:
V\in H^t(\mathbb{R}^n),\ V=v\text{ in }\Omega\right\}.
\end{equation}
Recall that, for a bounded Lipschitz domain $\Omega\subset\mathbb{R}^n$ and
$t\in(0,1)$, the norm in \eqref{e2.1} is equivalent to
\begin{equation}\label{e2.2}
\|v\|_{L^2(\Omega)}+\left[\iint_{\Omega\times \Omega}
\frac{|v(x)-v(y)|^2}{|x-y|^{n+2t}}\,dx\,dy\right]^{1/2}.
\end{equation}
For $t\in(0,1)$, let $[v]_{t,\Omega}$ denote the square root of the double
integral appearing in \eqref{e2.2}. Moreover, for $t\in(1,2)$, an equivalent
norm in \eqref{e2.1} is given by
$\|v\|_{H^1(\Omega)}+\sum_{i=1}^n\|\partial_iv\|_{H^{t-1}(\Omega)}$
(see, for instance, \cite{DNPV12,McLean00}).

Choose a radial function $\vartheta_0\in C_{\rm c}^\infty(\mathbb{R}^n)$
such that $0\le\vartheta_0\le1$, $\vartheta_0\equiv1$ on $B({\bf0},1)$, and
$\operatorname{supp}(\vartheta_0)\subset B({\bf0},2)$. For any $k\in\mathbb{N}$
and $\xi\in\mathbb{R}^n$, let
\begin{equation*}
\vartheta_k(\xi):=\vartheta_0(2^{-k}\xi)-\vartheta_0(2^{1-k}\xi).
\end{equation*}
Then $\{\vartheta_k\}_{k\in\mathbb{N}_0}$ forms an inhomogeneous smooth dyadic partition of unity. For any $k\in\mathbb{N}_0$, let
$\Delta_kf:=\mathcal F^{-1}(\vartheta_k\mathcal F f)$.
For any $f\in\mathcal S'(\mathbb{R}^n)$, $t\in\mathbb{R}$, and $p\in[1,\infty]$,
the Besov norm $\|f\|_{B^t_{2,p}(\mathbb{R}^n)}$ of $f$
on $\mathbb{R}^n$ is defined by setting
\begin{equation}\label{equ2.3}
\|f\|_{B^t_{2,p}(\mathbb{R}^n)}:=
\begin{cases}
\displaystyle\left(\sum_{k=0}^\infty
\left[2^{kt}\|\Delta_kf\|_{L^2(\mathbb{R}^n)}\right]^p\right)^{1/p},&p<\infty,\\[3pt]
\displaystyle\sup_{k\in\mathbb{N}_0}2^{kt}\|\Delta_kf\|_{L^2(\mathbb{R}^n)},&p=\infty.
\end{cases}
\end{equation}
The space $B^t_{2,p}(\Omega)$ is defined by restriction, with the
corresponding infimum norm as in \eqref{e2.1}. In particular,
$B^t_{2,2}(\Omega)=H^t(\Omega)$ with equivalent norms
(see, for instance, \cite{Tri02}).

Let $X_0$ and $X_1$ be two Banach spaces continuously embedded into a common
Hausdorff topological vector space. Such a pair $(X_0,X_1)$ is said to be a
\emph{compatible couple}. For any $f\in X_0+X_1$ and $\tau\in(0,\infty)$, the
\emph{$K$-functional} $\mathscr K(\tau,f,X_0,X_1)$ is defined by setting
\begin{equation*}
\mathscr K(\tau,f,X_0,X_1):=\inf_{{\genfrac{}{}{0pt}{}{f=f_0+f_1}
{f_0\in X_0,\ f_1\in X_1}}}
\left(\|f_0\|_{X_0}+\tau\|f_1\|_{X_1}\right).
\end{equation*}
For any $\theta\in(0,1)$ and $p\in[1,\infty]$, the \emph{real interpolation
norm} $\|f\|_{(X_0,X_1)_{\theta,p}}$ is defined by setting, when $p\in[1,\infty)$,
\begin{equation*}
\|f\|_{(X_0,X_1)_{\theta,p}}:=\left(\int_0^\infty
\left[\tau^{-\theta}\mathscr K(\tau,f,X_0,X_1)\right]^p\,\frac{d\tau}{\tau}\right)^{1/p},
\end{equation*}
and, when $p=\infty$,
\begin{equation*}
\|f\|_{(X_0,X_1)_{\theta,\infty}}:=\sup_{\tau\in(0,\infty)}\tau^{-\theta}
\mathscr K(\tau,f,X_0,X_1).
\end{equation*}
For any $\theta\in(0,1)$ and $p\in[1,\infty]$, the \emph{real interpolation
space} $(X_0,X_1)_{\theta,p}$ is defined to be the set of all
$f\in X_0+X_1$ satisfying $\|f\|_{(X_0,X_1)_{\theta,p}}<\infty$
equipped with the norm $\|\cdot\|_{(X_0,X_1)_{\theta,p}}$
(see, for instance, \cite{BL76}).

Then we have the following conclusions for Sobolev spaces and Besov spaces.

\begin{lemma}\label{l1}
Let $n\ge2$, and let $\Omega\subset\mathbb{R}^n$ be a bounded Lipschitz
domain. Then the following statements hold.
\begin{itemize}
\item[\rm(i)] There exists a linear extension operator $E$ such that, for any
$t\in\mathbb{R}$ and $p\in[1,\infty]$, the operators
\begin{equation*}
E:H^t(\Omega)\longrightarrow H^t(\mathbb{R}^n)
\quad\text{and}\quad
E:B^t_{2,p}(\Omega)\longrightarrow B^t_{2,p}(\mathbb{R}^n)
\end{equation*}
are bounded and $\mathcal RE=I$, where $\mathcal R$ denotes the restriction
to $\Omega$. Moreover, for any fixed $\chi\in C_{\rm c}^\infty(\mathbb{R}^n)$,
the multiplication operator $f\mapsto\chi f$ is bounded on these spaces.

\item[\rm(ii)] Let $t\in(0,1)$. A function $w\in H^t(\Omega)$ belongs to
$H^{t+1}(\Omega)$ if and only if each distributional derivative
$\partial_iw$, with $i\in\{1,\ldots,n\}$, belongs to $H^t(\Omega)$. Moreover,
\begin{equation}\label{e2.3}
\|w\|_{H^{t+1}(\Omega)}\sim\|w\|_{H^t(\Omega)}
+\sum_{i=1}^n\|\partial_iw\|_{H^t(\Omega)}
\end{equation}
with equivalent norms.

\item[\rm(iii)] Let $t\in(0,1/2)$. Then, for any $w\in H^t(\Omega)$,
\begin{equation}\label{e2.4}
\int_\Omega\delta^{-2t}|w|^2\,dx\lesssim\|w\|_{H^t(\Omega)}^2,
\end{equation}
where, for any $x\in\Omega$, $\delta(x):=\mathrm{dist\,}(x,\partial\Omega)$. Furthermore, for any
$w\in H^t(\mathbb{R}^n_+)$ with $t\in(0,1)$, if either $t\in(0,1/2)$ or $t\in(1/2,1)$ with
$\operatorname{Tr}w=0$, then
\begin{equation}\label{e2.5}
\int_{\mathbb{R}^n_+} \rho^{-2t}|w(x',\rho)|^2\,dx'\,dt
\lesssim\|w\|_{H^t(\mathbb{R}^n_+)}^2.
\end{equation}
Here $\operatorname{Tr}$ denotes the boundary trace operator, which agrees
with $w(x',0)$ for smooth $w$. For $t\in(1/2,1)$ and
$w\in H^t(\mathbb{R}^n_+)$, the trace operator $\operatorname{Tr}$ satisfies
\begin{equation}\label{e2.6}
\|\operatorname{Tr}w\|_{H^{t-1/2}(\mathbb{R}^{n-1})}\lesssim\|w\|_{H^t(\mathbb{R}^n_+)}.
\end{equation}

\item[\rm(iv)] For any $t\in\mathbb{R}$, $\theta\in(0,1)$, and $p\in[1,\infty]$,
\begin{equation}\label{e2.7}
\left(H^t(\Omega),H^{t+1}(\Omega)\right)_{\theta,p}=B^{t+\theta}_{2,p}(\Omega),
\end{equation}
and, for any $\tau\in(0,\infty)$ and $\theta\in(0,1)$,
\begin{equation}\label{e2.8}
\left(H^t(\Omega),B^{t+\tau}_{2,\infty}(\Omega)\right)_{\theta,2}=H^{t+\theta\tau}(\Omega).
\end{equation}
The identities \eqref{e2.7} and \eqref{e2.8} hold with equivalent norms.
\end{itemize}
\end{lemma}

\begin{proof}
We first show (i). By the common extension operator $E$ constructed in \cite{Ry99}
and \cite[(1.6), Theorem 1.5(1), and Lemma 6.3]{SY24}, we find that, for
$X:=H^t$ or $X:=B^t_{2,p}$, $\mathcal RE=I$ and
\begin{equation*}
\|\chi f\|_{X(\Omega)}\le\|\chi Ef\|_{X(\mathbb{R}^n)}
\lesssim\|Ef\|_{X(\mathbb{R}^n)}\lesssim\|f\|_{X(\Omega)}.
\end{equation*}
This finishes the proof of (i).

We now prove (ii). From the derivative norm theorem of Shi and Yao
\cite[Theorem 1.1, (1.2), and Proposition 6.4]{SY24} with smoothness order $t+1$,
it follows that
\begin{equation*}
\|w\|_{H^{t+1}(\Omega)}\sim\sum_{|\alpha|\le1}\|\partial^\alpha w\|_{H^t(\Omega)}.
\end{equation*}
This, together with the fact that $\partial^0w=w$, implies that \eqref{e2.3}
holds with equivalent norms. This shows (ii).

We next prove (iii). The estimate \eqref{e2.4} is justified in
\cite[proof of Proposition 3.8, Case 1]{GH24}. From \cite[Theorem 1.1]{FS10}
with $p=2$ and fractional order $t$, and from the density argument in
\cite[Sections 6.5, 6.6, and 9.3]{Le23}, we deduce that \eqref{e2.5} holds.
By the full-space trace estimate in \cite[Proposition 3.8]{DNPV12} and the
restriction norm \eqref{e2.1} with $\Omega=\mathbb{R}^n_+$, we find that
\eqref{e2.6} holds for $t\in(1/2,1)$. This finishes the proof of (iii).

We finally prove (iv). From the full-space interpolation theorem
\cite[Section 6.4]{BL76} and the common extension operator $E$ constructed in \cite{Ry99}
and \cite[(1.6), Theorem 1.5(1), and Lemma 6.3]{SY24},
we infer that \eqref{e2.7} and \eqref{e2.8} hold with equivalent norms.
This finishes the proof of (iv) and hence  Lemma \ref{l1}.
\end{proof}

\subsection{Energy estimates}\label{s2.2}

For any $x,y\in\mathbb{R}^n$ with $x\ne y$, let
$K(x,y):=|x-y|^{-n-2s}.$ For any $z\in\Omega^c$ and any function $v$ on $\Omega$, the function
$\rho$ is defined by setting
\begin{equation}\label{e2.9}
\rho(z):=\int_\Omega K(z,x)\,dx,
\end{equation}
and the operator $B$ is defined by setting
\begin{equation*}
Bv(z):=\rho(z)^{-1}\int_\Omega K(z,x)v(x)\,dx,
\end{equation*}
whenever the above integrals are well defined.

\begin{lemma}\label{l2}
Let $n\ge2$, let $s\in(0,1)$, and let $\Omega\subset\mathbb{R}^n$ be a bounded
Lipschitz domain. Then the following assertions hold.
\begin{itemize}
\item[\rm(i)] $\sigma_s(\Omega^c)<\infty$ and, for any $z\in\Omega^c$,
\begin{equation}\label{e2.10}
a_s(z)\le C\rho(z),
\end{equation}
where $a_s(z):=(1-s)\frac{\delta(z)^{-s}}{1+\delta(z)^{n+s}}$.
Moreover, there exist constants $r_*\in(0,r_0)$ and
$R_*\in(2\mathrm{diam}\,(\Omega),\infty)$ such that, for any $z\in\Omega^c$ with
$\delta(z)\in(0,r_*)$,
\begin{equation}\label{e2.11}
C^{-1}\delta(z)^{-2s}\le\rho(z)\le C\delta(z)^{-2s},
\end{equation}
and, for any $x\in\Omega$ and $z\in\Omega^c$ with $|z-x|>R_*$,
\begin{equation}\label{e2.12}
C^{-1}|z-x|^{-n-2s}\le\rho(z)\le C|z-x|^{-n-2s}.
\end{equation}
Here the positive constants $C$, $r_*$, and $R_*$ depend only on $n$, $s$,
and $\Omega$.

\item[\rm(ii)] For any $x\in\Omega$,
\begin{equation}\label{e2.13}
J(x):=\int_{\Omega^c}\frac{a_s(z)}{\rho(z)}K(z,x)\,dz
\le C[1+\delta(x)^{-s}],
\end{equation}
where the positive constant $C$ depends only on $n$, $s$, and $\Omega$.

\item[\rm(iii)] The operator $B$ maps $H^{s/2}(\Omega)$ continuously
into $L^2(\Omega^c,\sigma_s)$. Precisely, for any $w\in H^{s/2}(\Omega)$,
\begin{equation}\label{e2.14}
\|Bw\|_{L^2(\Omega^c,\sigma_s)}\le C\|w\|_{H^{s/2}(\Omega)},
\end{equation}
where the positive constant $C$ depends only on $n$, $s$, and $\Omega$.
\end{itemize}
\end{lemma}

\begin{proof}
We first prove (i).
By a finite graph cover near the boundary $\partial\Omega$ and the fact that,
for any $x\in\Omega$ and $z\in\Omega^c$ satisfying that $|z-x|\ge2\mathrm{diam\,}(\Omega)$,
$\delta(z)\sim |z-x|$, we conclude that
\begin{align*}
\sigma_s(\Omega^c)&\leq C\int_0^{r_0}r^{-s}\,dr+C
+C\int_{2\mathrm{diam\,}(\Omega)}^\infty r^{-1-2s}\,dr\\
&=Cr_0^{1-s}+C+C\mathrm{diam\,}(\Omega)^{-2s}<\infty,
\end{align*}
where $r_0\in(0,\infty)$ is as in Definition \ref{d2.1}.

Let $z\in\Omega^c$ and $r_*\in(0,\,r_0)$. From the fact that, for any $x\in\Omega$, $|z-x|\geq \delta(z)$, we deduce that
\begin{equation}\label{e2.15}
\rho(z)\leq\int_{\{y\in\mathbb{R}^n:|y-z|\geq \delta(z)\}}|z-y|^{-n-2s}\,dy
\lesssim\delta(z)^{-2s}.
\end{equation}
Meanwhile, by the interior cone condition of $\Omega$, we find that there
exists a constant $c\in(0,1)$ such that, for any $z\in\Omega^c$ with
$0<\delta(z)<r_*$ and a nearest boundary point $p$, there is a ball
$B(x_p,c\delta(z))\subset\Omega$ with $|x_p-p|\le\delta(z)$. From the fact
that, for any $x\in B(x_p,c\delta(z))$, $|z-x|\le C\delta(z)$, it follows that
\begin{equation*}
\rho(z)\ge\int_{B(x_p,c\delta(z))}K(z,x)\,dx
\ge C\delta(z)^n\delta(z)^{-n-2s}=C\delta(z)^{-2s}.
\end{equation*}
By this and \eqref{e2.15}, we conclude that \eqref{e2.11} holds.

Fix $x_0\in\Omega$ and $R_*\in(2\,\mathrm{diam\,}(\Omega),\infty)$.
From the fact that, for any $x\in\Omega$ and $z\in\Omega^c$ satisfying
$|z-x_0|>R_*$, $|z-x|\sim|z-x_0|$, we deduce that \eqref{e2.12}
holds. On the remaining exterior set, $\delta(z)\ge r_*$ and
$|z-x_0|\le R_*$, the function $a_s$ is bounded and $\rho$ is bounded
below by a positive constant. By this, \eqref{e2.11}, and \eqref{e2.12},
we find that \eqref{e2.10} holds. This finishes the proof of (i).

Next, we prove (ii). By (i), we have $a_s(z)/\rho(z)\le C\delta(z)^s$ near
the boundary; also, for any $x\in\Omega$ and $z\in\Omega^c$,
$\max\{\delta(z),\delta(x)\}\le|z-x|$. From these facts, it follows that
\begin{align}\label{e2.16}
\int_{\{z\in\Omega^c:\delta(z)\le r_*\}}\frac{a_s(z)}{\rho(z)}K(z,x)\,dz
&\le C\int_{\{z\in\mathbb{R}^n:|z-x|\ge\delta(x)\}}|z-x|^{-n-s}\,dz
\nonumber\\
&=C\int_{\delta(x)}^\infty r^{-1-s}\,dr\le C\delta(x)^{-s}.
\end{align}
By \eqref{e2.10} and \eqref{e2.12}, we conclude that, for any $x\in\Omega$,
\begin{equation}\label{e2.17}
\int_{\{z\in\Omega^c:\delta(z)\ge r_*\}}\frac{a_s(z)}{\rho(z)}K(z,x)\,dz
\le C+C\int_{\{z\in\Omega^c:|z-x|>R_*\}}
|z-x|^{-n-2s}\,dz\le C.
\end{equation}
From \eqref{e2.16} and \eqref{e2.17}, we deduce that \eqref{e2.13} holds.
This finishes the proof of (ii).

Finally, we prove (iii). By H\"older's inequality, we find that, for any
$z\in\Omega^c$,
\begin{equation*}
|Bw(z)|^2\le\rho(z)^{-1}\int_\Omega K(z,x)|w(x)|^2\,dx.
\end{equation*}
From this and \eqref{e2.13}, it follows that
\begin{equation*}
\|Bw\|_{L^2(\Omega^c,\sigma_s)}^2
\le\int_\Omega|w(x)|^2J(x)\,dx
\le C\int_\Omega\left[1+\delta^{-s}\right]|w|^2\,dx.
\end{equation*}
By Hardy's inequality \eqref{e2.4} with $t:=s/2<1/2$, we conclude that
\begin{equation*}
\|Bw\|_{L^2(\Omega^c,\sigma_s)}^2\le C\left[\|w\|_{L^2(\Omega)}^2
+\int_\Omega\delta^{-s}|w|^2\,dx\right]
\le C\|w\|_{H^{s/2}(\Omega)}^2,
\end{equation*}
which proves \eqref{e2.14}. This finishes the proof of (iii) and hence Lemma \ref{l2}.
\end{proof}

\begin{lemma}\label{l3}
Let $n\ge2$, let $s\in(0,1)$, and let $\Omega\subset\mathbb{R}^n$ be a bounded
Lipschitz domain. Assume that $\mathcal{V}$ and $\mathcal{V}_g$ are the
spaces defined in \eqref{e1.6} and \eqref{e1.8}, respectively. Then the
following assertions hold.
\begin{itemize}
\item[\rm(i)] The space $\mathcal{V}$ is a Hilbert space, and, for any
$v\in\mathcal{V}$,
\begin{equation}\label{e2.18}
\|v\|_{H^s(\Omega)}+\|v\|_{L^2(\Omega^c,\sigma_s)}
\le C\|v\|_{\mathcal{V}},
\end{equation}
where the positive constant $C$ depends only on $n$, $s$, and $\Omega$.
Moreover, if $v$ belongs to the closed subspace
$\widetilde{\mathcal{V}}:=\{v\in\mathcal{V}:\int_\Omega v\,dx=0\},$
then the Poincar\'e type inequality
\begin{equation*}
\|v\|_{H^s(\Omega)}+\|v\|_{L^2(\Omega^c,\sigma_s)}
\le C\mathcal{E}_\Omega(v,v)^{1/2}
\end{equation*}
holds, with the positive constant $C$ depending only on $n$, $s$, and
$\Omega$.

\item[\rm(ii)] For any $g\in\mathcal X_s(\Omega^c)$ satisfying
$\int_{\Omega^c}g\,dz=0$, there exists a unique weak solution
$u\in\widetilde{\mathcal{V}}$ to problem \eqref{e1.1} in the sense of
Definition~\ref{def:weak}, and
\begin{equation}\label{e2.19}
\mathcal{E}_\Omega(u,u)^{1/2}+\|u\|_{H^s(\Omega)}
+\|u\|_{L^2(\Omega^c,\sigma_s)}
\le C\|g\|_{\mathcal X_s(\Omega^c)},
\end{equation}
where the positive constant $C$ depends only on $n$, $s$, and $\Omega$.

\item[\rm(iii)] For any function $g$ as in \textup{(ii)}, the space
$\mathcal{V}_g$ is dense in $\mathcal{V}$. Moreover, for $u\in
\widetilde{\mathcal{V}}$, \eqref{e1.11} holds for any $v\in\mathcal{V}_g$
if and only if it holds for any $v\in\mathcal{V}$.
\end{itemize}
\end{lemma}

\begin{proof}
We first prove (i). The Hilbert space structure and the Poincar\'e type
inequality are standard (see, for instance, \cite[Proposition 3.1 and
Lemma 3.10]{DRV17}). The exterior $L^2(\sigma_s)$ bound in \eqref{e2.18} was
also proved in \cite[Remark 1.3 and the proof of Theorem 1.2(1)]{GH24},
with the equivalent density $a_s$ in place of $\tau_s$. This finishes the
proof of (i).

Next, we prove (ii). By H\"older's inequality and \eqref{e2.18}, we find
that, for $g=a_sq$ and $v\in\mathcal{V}$,
\begin{equation}\label{e2.20}
\left|\int_{\Omega^c}gv\,dz\right|
\le\|q\|_{L^2(\Omega^c,\sigma_s)}
\|v\|_{L^2(\Omega^c,\sigma_s)}
\le C\|q\|_{L^2(\Omega^c,\sigma_s)}\|v\|_{\mathcal{V}}.
\end{equation}
From (i), we deduce that, for any $u,v\in\widetilde{\mathcal{V}}$,
\begin{equation*}
|\mathcal{E}_\Omega(u,v)|\le\|u\|_{\mathcal{V}}\|v\|_{\mathcal{V}}
\quad\text{and}\quad
\mathcal{E}_\Omega(u,u)\ge c\|u\|_{\mathcal{V}}^2.
\end{equation*}
This, together with the Lax--Milgram theorem, implies that there exists a
unique $u\in\widetilde{\mathcal{V}}$ such that, for any
$w\in\widetilde{\mathcal{V}}$,
$\mathcal{E}_\Omega(u,w)=\int_{\Omega^c}gw\,dz.$
For any $v\in\mathcal{V}$, let
$v_\Omega:=|\Omega|^{-1}\int_\Omega v\,dx$. Since
$v-v_\Omega\in\widetilde{\mathcal{V}}$ and
$\int_{\Omega^c}g\,dz=0$, it follows that
\begin{equation*}
\mathcal{E}_\Omega(u,v)=\mathcal{E}_\Omega(u,v-v_\Omega)
=\int_{\Omega^c}g\,[v-v_\Omega]\,dz=\int_{\Omega^c}gv\,dz.
\end{equation*}
This proves the weak identity and hence the existence in (ii). The
uniqueness follows from the same coercive estimate
$\mathcal{E}_\Omega(u,u)\ge c\|u\|_{\mathcal{V}}^2$. By \eqref{e2.20} and
(i), we find that
\begin{equation*}
\mathcal{E}_\Omega(u,u)\le C\|q\|_{L^2(\Omega^c,\sigma_s)}\|u\|_{\mathcal{V}}
\le C\|q\|_{L^2(\Omega^c,\sigma_s)}\mathcal{E}_\Omega(u,u)^{1/2}.
\end{equation*}
From this, (i), and \eqref{e1.4}, we infer that \eqref{e2.19} holds. This
finishes the proof of (ii).

Finally, we prove (iii). By H\"older's inequality and the fact that
$\sigma_s(\Omega^c)<\infty$, we conclude that
\begin{equation*}
\|g\|_{L^1(\Omega^c)}\le\sigma_s(\Omega^c)^{1/2}\|q\|_{L^2(\Omega^c,\sigma_s)}<\infty.
\end{equation*}
For any $v\in\mathcal{V}$, define its bounded truncation by
$T_kv:=\max\{-k,\min\{v,k\}\}$ with $k\in\mathbb{N}$. From the contraction
property, we deduce that $|T_kv|\le k$,
$\|T_kv\|_{L^2(\Omega)}\le\|v\|_{L^2(\Omega)}$,
$\mathcal{E}_\Omega(T_kv,T_kv)\le\mathcal{E}_\Omega(v,v)$, and
\begin{equation*}
\int_{\Omega^c}|g|\,|T_kv|^2\,dz
\le k^2\|g\|_{L^1(\Omega^c)}<\infty.
\end{equation*}
Thus $T_kv\in\mathcal{V}_g$. By the same contraction property, we conclude
that, for any $x,y\in\mathbb{R}^n$,
\begin{equation*}
\left|[v(x)-T_kv(x)]-[v(y)-T_kv(y)]\right|\le|v(x)-v(y)|.
\end{equation*}
From this and the dominated convergence theorem, we infer that
$\|v-T_kv\|_{\mathcal{V}}\to0$ as $k\to\infty$, which, together with
\eqref{e2.18}, further implies that
\begin{equation}\label{e2.21}
\int_{\Omega^c}|g(z)|\,|v(z)-T_kv(z)|\,dz
\le C\|q\|_{L^2(\Omega^c,\sigma_s)}
\|v-T_kv\|_{\mathcal{V}}\to0
\end{equation}
as $k\to\infty$. By H\"older's inequality, we conclude that
\begin{equation}\label{e2.22}
|\mathcal{E}_\Omega(u,v-T_kv)|\le\mathcal{E}_\Omega(u,u)^{1/2}
\mathcal{E}_\Omega(v-T_kv,v-T_kv)^{1/2}\to0
\end{equation}
as $k\to\infty$. From \eqref{e2.21} and \eqref{e2.22}, it follows that, for
any $v\in\mathcal{V}$,
\begin{equation*}
\mathcal{E}_\Omega(u,v)=\lim_{k\to\infty}\mathcal{E}_\Omega(u,T_kv)
=\lim_{k\to\infty}\int_{\Omega^c}gT_kv\,dz=\int_{\Omega^c}gv\,dz.
\end{equation*}
This finishes the proof of (iii) and hence Lemma \ref{l3}.
\end{proof}

For any $v\in H^s(\Omega)$, its \emph{exterior extension} $\mathcal Mv$ is
defined by setting
\begin{equation*}
(\mathcal Mv)(x):=
\begin{cases}
v(x),&x\in\Omega,\\
Bv(x),&x\in\Omega^c.
\end{cases}
\end{equation*}
The same definition applies to any function $u$ already defined on
$\mathbb{R}^n$, by using only its values in $\Omega$. In this case,
$\mathcal Mu$ is a new global function which equals $u$ in $\Omega$ and
equals $Bu$ on $\Omega^c$, and hence differs from $u$ outside $\Omega$.
For any $v,w\in H^s(\Omega)$, the bilinear form $\mathcal E_N$ is defined by
setting
\begin{equation*}
\mathcal E_N(v,w):=\mathcal E_\Omega(\mathcal Mv,\mathcal Mw),
\end{equation*}
and the bilinear form $\mathfrak a$ is defined by setting
\begin{equation*}
\mathfrak a(v,w):=\mathcal E_N(v,w)+(v,w)_{L^2(\Omega)}.
\end{equation*}

Inspired by \cite[Lemma A.2 and Proposition A.3]{AFR23},
we have the following lemma.

\begin{lemma}\label{l4}
Let $n\ge2$, let $s\in(0,1)$, and let $\Omega\subset\mathbb{R}^n$ be a bounded
Lipschitz domain. Then the following assertions hold.
\begin{itemize}
\item[\rm(i)] For any $v\in H^s(\Omega)$, $\mathcal Mv\in\mathcal{V}$ is the
unique minimizer of $\mathcal E_\Omega(u,u)$ among all $u\in\mathcal{V}$
that coincide with $v$ on $\Omega$. Moreover,
\begin{equation}\label{e2.23}
\mathfrak a(v,v)\sim\|v\|_{H^s(\Omega)}^2
\end{equation}
with equivalent norms.

\item[\rm(ii)] For any $v\in\mathcal{V}$, $v-Bv\in
L^2(\Omega^c,\rho\,dz)$, and, for any $u,v\in\mathcal{V}$,
\begin{equation}\label{e2.24}
\mathcal E_\Omega(u,v)=\mathcal E_N(u,v)+c_{(n,s)}\int_{\Omega^c}\rho
(u-Bu)(v-Bv)\,dz.
\end{equation}

\item[\rm(iii)] Let $q\in L^2(\Omega^c,\sigma_s)$ satisfy
$\int_{\Omega^c}q\,d\sigma_s=0$, and let $g:=a_sq$. Assume that
$u:\mathbb{R}^n\to\mathbb{R}$ is measurable and that
$u\in H^s(\Omega)$. Then $u$ is a weak solution in the sense of
Definition~\ref{def:weak} if and only if $\int_\Omega u\,dx=0$,
for any $z\in\Omega^c$,
\begin{equation}\label{e2.25}
u(z)=Bu(z)+\frac{a_s(z)q(z)}{c_{(n,s)}\rho(z)}
\end{equation}
and, for any $w\in H^s(\Omega)$,
\begin{equation}\label{e2.26}
 \mathcal{E}_N(u,w)=F_q(w):=\int_{\Omega^c}q(z)Bw(z)\,d\sigma_s(z).
\end{equation}

\item[\rm(iv)] For any $q$ as in \textup{(iii)}, the functional $F_q$
defined in \eqref{e2.26} is continuous on $H^{s/2}(\Omega)$, and
\begin{equation}\label{e2.27}
\|F_q\|_{(H^{s/2}(\Omega))^*}
\le C\|q\|_{L^2(\Omega^c,\sigma_s)},
\end{equation}
where the positive constant $C$ depends only on $n$, $s$, and $\Omega$.
\end{itemize}
\end{lemma}

\begin{proof}
We first prove (i).
Let $v\in H^s(\Omega)$ and $u:\mathbb{R}^n\to\mathbb{R}$ be measurable
with $u=v$ in $\Omega$.
From the definition of $B$, we deduce that, for any
$z\in\Omega^c$,
\begin{equation}\label{e2.28}
 \int_\Omega K(z,x)[u(x)-Bu(z)]\,dx=0.
\end{equation}
By this and the definition
of $B$, we find that, for any $z\in\Omega^c$,
\begin{equation}\label{e2.29}
\int_\Omega K(z,x)|u(z)-u(x)|^2\,dx=\rho(z)|u(z)-Bu(z)|^2
+\int_\Omega K(z,x)|u(x)-Bu(z)|^2\,dx.
\end{equation}
From \eqref{e2.29}, it follows that, for any function
$u$ equal to $v$ in $\Omega$,
\begin{equation}\label{e2.30}
\mathcal{E}_\Omega(u,u)=\mathcal{E}_N(v,v)+c_{(n,s)}\int_{\Omega^c}
\rho(z)|u(z)-Bv(z)|^2\,dz.
\end{equation}
By Lemma~\ref{l1}\textup{(i)} and \eqref{e2.30}, we conclude that
\begin{equation*}
\mathcal{E}_N(v,v)\leq\mathcal{E}_\Omega(Ev,Ev)\leq C\|v\|_{H^s(\Omega)}^2
\end{equation*}
and $\mathcal M v\in\mathcal{V}$.
From \eqref{e2.30}, we deduce that, for any $u\in\mathcal{V}$ equal to $v$
in $\Omega$, $\mathcal{E}_\Omega(u,u)\geq\mathcal{E}_N(v,v)$, with equality
if and only if $u=Bv$ almost everywhere on $\Omega^c$.

By Lemma~\ref{l1}\textup{(i)}, \eqref{e2.2}, and the minimum
energy estimate of $\mathcal{E}_\Omega(u,u)$ with $u=v$ in $\Omega$ proved above, we find that,
for any $v\in H^s(\Omega)$ and its extension $Ev\in H^s(\mathbb{R}^n)$,
\begin{equation*}
c\|v\|_{H^s(\Omega)}^2\leq\|v\|_{L^2(\Omega)}^2
+\frac{c_{(n,s)}}2\iint_{\Omega\times\Omega}K(x,y)|v(x)-v(y)|^2\,dx\,dy
\leq\mathfrak a(v,v)
\end{equation*}
and
\begin{equation*}
\mathfrak a(v,v)\leq\|v\|_{L^2(\Omega)}^2+\mathcal{E}_\Omega(Ev,Ev)
\leq C\|Ev\|_{H^s(\mathbb{R}^n)}^2\leq C\|v\|_{H^s(\Omega)}^2.
\end{equation*}
From this, it follows that \eqref{e2.23} holds, which proves
\textup{(i)}.

Now, we show (ii). By \eqref{e2.18}, \eqref{e2.30}, and
\textup{(i)}, we conclude that, for any $v\in\mathcal{V}$,
\begin{align*}
c_{(n,s)}\int_{\Omega^c}\rho(z)|v(z)-Bv(z)|^2\,dz
=\mathcal{E}_\Omega(v,v)-\mathcal{E}_N(v,v)\leq\mathcal{E}_\Omega(v,v)<\infty.
\end{align*}
This proves the weighted integrability asserted in \textup{(ii)}. For any
function $k\in L^2(\Omega^c,\rho\,dz)$ on $\Omega^c$, let $Z_k$ be equal to zero in
$\Omega$ and to $k$ on $\Omega^c$.
From the definition of $\mathcal{E}_\Omega$, we deduce that
\begin{equation}\label{e2.31a}
\mathcal{E}_\Omega(Z_k,Z_k)=c_{(n,s)}\|k\|_{L^2(\Omega^c,\rho\,dz)}^2
\quad\text{and}\quad
Z_k\in\mathcal{V}.
\end{equation}
For any $v\in\mathcal{V}$, take
$k:=v-Bv$ on $\Omega^c$. Then $v=\mathcal M v+Z_k$.
By \textup{(i)} and \eqref{e2.30}, we find that
$$
\|v\|_{\mathcal{V}}^2 =\mathfrak a(v,v) +c_{(n,s)}\|k\|_{L^2(\Omega^c,\rho\,dz)}^2.
$$
For any $u,v\in\mathcal{V}$, let $r_u:=u-Bu$ and $r_v:=v-Bv$ on
$\Omega^c$. Then from this, it follows that, for any
pair $(z,x)\in\Omega^c\times\Omega$,
\begin{equation}\label{e2.32}
[u(z)-u(x)][v(z)-v(x)]=[r_u(z)+Bu(z)-u(x)][r_v(z)+Bv(z)-v(x)].
\end{equation}
By \eqref{e2.28}, we conclude that, for any $z\in\Omega^c$,
\begin{equation}\label{e2.33}
r_u(z)\int_\Omega K(z,x)[Bv(z)-v(x)]\,dx=0
\end{equation}
and
\begin{equation}\label{e2.34}
r_v(z)\int_\Omega K(z,x)[Bu(z)-u(x)]\,dx=0.
\end{equation}
From \eqref{e2.32}, \eqref{e2.33}, and \eqref{e2.34}, we deduce that, for
any $z\in\Omega^c$,
\begin{align}\label{e2.35}
&\int_\Omega K(z,x)[u(z)-u(x)][v(z)-v(x)]\,dx\nonumber\\
&\quad=\rho(z)r_u(z)r_v(z)
+\int_\Omega K(z,x)[Bu(z)-u(x)][Bv(z)-v(x)]\,dx.
\end{align}
Integrating \eqref{e2.35} over $\Omega^c$ and adding the resulting identity
to the integral over $\Omega\times\Omega$, we conclude that \eqref{e2.24}
holds. This proves (ii).

We next prove (iii). For any $z\in\Omega^c$, let
$b_g(z):=\frac{g(z)}{c_{(n,s)}\rho(z)}$.
By $g=a_sq$ and \eqref{e2.10}, we find that
\begin{align}\label{e2.39}
\|b_g\|_{L^2(\Omega^c,\rho\,dz)}^2=c_{(n,s)}^{-2}\int_{\Omega^c}\frac{|g|^2}{\rho}\,dz
=c_{(n,s)}^{-2}\int_{\Omega^c}\frac{a_s}{\rho}|q|^2\,
d\sigma_s\leq C\|q\|_{L^2(\Omega^c,\sigma_s)}^2.
\end{align}
From H\"older's inequality and \eqref{e2.39}, it follows that, for any
$k\in L^2(\Omega^c,\rho\,dz)$,
\begin{equation}\label{e2.40}
\int_{\Omega^c}|gk|\,dz\leq c_{(n,s)}\|b_g\|_{L^2(\Omega^c,\rho\,dz)}
\|k\|_{L^2(\Omega^c,\rho\,dz)}.
\end{equation}

We first prove necessity. Suppose that $u$ is a weak solution in the
sense of Definition~\ref{def:weak}, so that
$u\in\widetilde{\mathcal{V}}$. By Lemma~\ref{l3}\textup{(iii)}, \eqref{e2.31a}, and
\eqref{e2.24}, we conclude that, for any
$k\in L^2(\Omega^c,\rho\,dz)$,
\begin{equation}\label{e2.41}
c_{(n,s)}\int_{\Omega^c}\rho[u-Bu]k\,dz
=\int_{\Omega^c}gk\,dz.
\end{equation}
From \textup{(ii)}, \eqref{e2.39}, and \eqref{e2.41} with taking
$k=u-Bu-b_g$, we deduce that
$k\in L^2(\Omega^c,\rho\,dz)$ and
$c_{(n,s)}\int_{\Omega^c}\rho|k|^2\,dz=0$.
By the positivity of $\rho$, we find that,
$u-Bu=b_g$,
which proves \eqref{e2.25}. In particular, $u=\mathcal M u$ when $g=0$.
From (i), Lemma~\ref{l3}\textup{(iii)}, and \eqref{e2.24}, it follows that, for any $w\in H^s(\Omega)$,
\begin{equation*}
\mathcal{E}_N(u,w)=\int_{\Omega^c}gBw\,dz
=\int_{\Omega^c}qBw\,d\sigma_s=F_q(w).
\end{equation*}
This proves \eqref{e2.26}.

We next prove sufficiency. Assume that $u\in H^s(\Omega)$,
$\int_\Omega u(x)dx=0$, and \eqref{e2.25} and \eqref{e2.26} hold.
By \textup{(i)}, \eqref{e2.25}, \eqref{e2.31a}, and \eqref{e2.39}, we conclude that
$u=\mathcal M u+Z_{b_g}\in\widetilde{\mathcal{V}}$.
For any $v\in\mathcal{V}$, let $k=v-Bv$ on $\Omega^c$.
From \textup{(ii)} and \eqref{e2.26}, we deduce that
\begin{align}\label{e2.42}
\mathcal{E}_\Omega(u,v)&=\mathcal{E}_N(u,v)+c_{(n,s)}\int_{\Omega^c}\rho b_gk\,dz
=F_q(v)+\int_{\Omega^c}gk\,dz\nonumber\\
&=\int_{\Omega^c}g[Bv+k]\,dz=\int_{\Omega^c}gv\,dz.
\end{align}
By \eqref{e2.40}, \eqref{e2.14}, and \eqref{e2.18}, we find that
every term in \eqref{e2.42} is absolutely convergent and that $u$ is a
weak solution in the sense of Definition~\ref{def:weak}.
This proves \textup{(iii)}.

We finally show (iv). By H\"older inequality and \eqref{e2.14}, we conclude that,
for any $w\in H^{s/2}(\Omega)$,
\begin{align}\label{e2.43}
|F_q(w)|\leq\|q\|_{L^2(\Omega^c,\sigma_s)}\|Bw\|_{L^2(\Omega^c,\sigma_s)}
\leq C\|q\|_{L^2(\Omega^c,\sigma_s)}\|w\|_{H^{s/2}(\Omega)}.
\end{align}
From the linearity of $B$ and \eqref{e2.43}, we deduce that
$F_q\in(H^{s/2}(\Omega))^*$ and that \eqref{e2.27} holds.
This finishes the proof of (iv) and hence Lemma~\ref{l4}.
\end{proof}

\section{Improved Sobolev regularity}\label{s3}

In this section, we first prove an energy comparison and
then establish the lifted Sobolev regularity estimate for the weak solution to \eqref{e1.1}.

\subsection{Energy comparison}\label{s3.1}

Let $n\geq2$, $s\in(0,1)$, and $\Omega\subset\mathbb{R}^n$ be a bounded Lipschitz domain.
By Definition \ref{d2.1}, we find that, for any boundary chart $\mathcal U$ satisfying that
$\mathcal U\subset B(x,r)$ with some $x\in\partial\Omega$ and $r\in(0,r_0)$,
$\Omega\cap\mathcal U=\{x\in\mathcal U:x_n>\varphi(x')\}$, where $\varphi$ is a Lipschitz function
satisfying $\|\nabla\varphi\|_{L^\infty(\mathbb{R}^{n-1})}\leq L$. Choose
$\eta\in C_{\rm c}^\infty(\mathcal U)$
with $0\leq\eta\leq1$, equal to one on a smaller neighborhood.
Define the inward cone $\mathcal C$ by
\begin{equation*}
\mathcal C:=\{h=(h',h_n)\in\mathbb{R}^n:h_n>(L+1)|h'|\}.
\end{equation*}
For $h\in\mathcal C$ with $0<|h|<h_0$, let
\begin{equation}\label{e3.2}
\Phi_h(x):=x+\eta(x)h
\end{equation}
and $T_hv:=v\circ\Phi_h.$ By the definition of $\Phi_h$, we conclude that,
for any $x\notin\mathrm{supp\,}\eta$, $\Phi_h(x)=x$.
Choose $h_0\leq1$ so that $q_h:=|h|\,\|\nabla\eta\|_{L^\infty(\mathbb{R}^n)}\leq1/2$.
From the Lipschitz bound for $\eta$, we deduce that, for any
$x,y\in\mathbb{R}^n$,
\begin{equation}\label{e3.4}
 (1-q_h)|x-y|\leq|\Phi_h(x)-\Phi_h(y)|\leq(1+q_h)|x-y|.
\end{equation}
By \eqref{e3.4}, we find that $\Phi_h$ is injective.
From the definition of $q_h$ and the contraction principle,
it follows that, for any $A\in\mathbb{R}^n$,
the map $x\mapsto A-\eta(x)h$ has a unique fixed point
$\Psi_h(A)$, which satisfies that
\begin{equation*}
\Psi_h(A)=A-\eta(\Psi_h(A))h,\quad
\Phi_h(\Psi_h(A))=A,\quad
\text{and}\quad
\Psi_h(\Phi_h(A))=A.
\end{equation*}
By $\eta\in C^\infty_{\rm c}(\mathcal U)$ and $q_h\leq1/2$, we conclude that, for any
$x,\xi\in\mathbb{R}^n$,
\begin{subequations}
\begin{equation}\label{e3.5a}
 D\Phi_h(x)\xi=\xi+h(\nabla\eta(x)\cdot\xi),
\end{equation}
\begin{equation}\label{e3.5b}
 |D\Phi_h(x)\xi|\geq(1-q_h)|\xi|,
\end{equation}
\begin{equation}\label{e3.5c}
 \ker D\Phi_h(x)=\{0\},
\end{equation}
\begin{equation}\label{e3.5d}
 \det D\Phi_h(x)\ne0,
\end{equation}
and
\begin{equation}\label{e3.5e}
 \|D\Phi_h-I\|_{L^\infty(\mathbb{R}^n)}\leq q_h\leq C|h|.
\end{equation}
\end{subequations}
From \eqref{e3.5a}--\eqref{e3.5e}, the smooth inverse function theorem, and global
bijectivity, we deduce that every local inverse is $\Psi_h$. By this, we find that, for any $A\in\mathbb{R}^n$,
$\Psi_h\in C^\infty(\mathbb{R}^n)$ and
$D\Psi_h(A)=[D\Phi_h(\Psi_h(A))]^{-1}$.
From the inward cone condition of $\Omega$, it follows that
\begin{equation*}
\Omega_h:=\Phi_h(\Omega)\subset\Omega.
\end{equation*}

The following lemma is an adaptation of the localized variational comparison
method of Savar\'e \cite{Sa98}. The smooth deformation and
change-of-variables framework are given in
\cite[Section 2.3 and Proposition 3.2]{BLN23} for zero exterior
Dirichlet data.

\begin{lemma}\label{l5}
Let $n\ge2$, let $s\in(0,1)$, and let $\Omega\subset\mathbb{R}^n$ be a bounded
Lipschitz domain. Then there exist positive constants $h_0$ and $C$, with
$C$ depending only on $n$, $s$, $\Omega$, and $\|\eta\|_{C^1(\mathbb{R}^n)}$,
such that, for any $v\in H^s(\Omega)$ and any $h\in\mathcal C$ satisfying
$0<|h|<h_0$, the function $T_hv:=v\circ\Phi_h$ belongs to $H^s(\Omega)$ and
\begin{equation*}
\mathfrak a(T_hv,T_hv)\le(1+C|h|)\,\mathfrak a(v,v).
\end{equation*}
\end{lemma}

\begin{proof}
For any open set $G$, let $\mathcal Q_G$ and $\mathcal E_G$ be defined as in
\eqref{e1.5} and \eqref{e1.7}, respectively, with $\Omega$ replaced by $G$.
If $G_1\subset G_2$, then, for any function $u$,
\begin{equation}\label{e3.8}
\mathcal E_{G_1}(u,u)\le\mathcal E_{G_2}(u,u).
\end{equation}
For any given $x,y\in\mathbb{R}^n$ and $h\in\mathcal C$ satisfying
$0<|h|<h_0$, let $A:=\Phi_h(x)$, $B:=\Phi_h(y)$, $\Psi_h:=\Phi_h^{-1}$, and
$J_h:=|\det D\Psi_h|$. From \eqref{e3.4} and \eqref{e3.5a}--\eqref{e3.5e},
we deduce that $(x,y)\in\mathcal Q_\Omega$ if and only if
$(A,B)\in\mathcal Q_{\Omega_h}$, and that
\begin{subequations}
\begin{equation}\label{e3.9a}
\Phi_h(\Omega^c)=\Omega_h^c,
\end{equation}
\begin{equation}\label{e3.9b}
|J_h-1|\le C|h|,
\end{equation}
\begin{equation}\label{e3.9c}
|\Psi_h(A)-\Psi_h(B)|\ge(1-C|h|)|A-B|
\end{equation}
for any $A,B\in\mathbb{R}^n$, and
\begin{equation}\label{e3.9d}
K_h(A,B):=\frac{J_h(A)J_h(B)}{|\Psi_h(A)-\Psi_h(B)|^{n+2s}}
\le(1+C|h|)K(A,B).
\end{equation}
\end{subequations}
By the change of variables, we obtain
\begin{equation}\label{e3.10}
\mathcal E_\Omega((\mathcal Mv)\circ\Phi_h,(\mathcal Mv)\circ\Phi_h)
=\frac{c_{(n,s)}}2\iint_{\mathcal Q_{\Omega_h}}
K_h(A,B)\big|(\mathcal Mv)(A)-(\mathcal Mv)(B)\big|^2\,dA\,dB.
\end{equation}
From \eqref{e3.10}, \eqref{e3.9a}--\eqref{e3.9d}, \eqref{e3.8}, and the
fact that $\Omega_h\subset\Omega$, it follows that, for
$W_h:=(\mathcal Mv)\circ\Phi_h$,
\begin{equation}\label{e3.11}
\mathcal E_\Omega(W_h,W_h)
\le(1+C|h|)\mathcal E_{\Omega_h}(\mathcal Mv,\mathcal Mv)
\le(1+C|h|)\mathcal E_N(v,v)<\infty.
\end{equation}
By the change of variables $A=\Phi_h(x)$, we conclude that
\begin{equation}\label{e3.12}
\|W_h\|_{L^2(\Omega)}^2
\le(1+C|h|)\int_{\Omega_h}|v|^2\,dx
\le(1+C|h|)\|v\|_{L^2(\Omega)}^2.
\end{equation}
From \eqref{e3.11}, \eqref{e3.12}, and Lemma~\ref{l3}\textup{(i)}, we
deduce that $W_h\in\mathcal V$, $T_hv\in H^s(\Omega)$, and $W_h=T_hv$ in
$\Omega$. This, together with Lemma~\ref{l4}(i), implies that
\begin{equation}\label{e3.13}
\begin{aligned}
\mathcal E_N(T_hv,T_hv)
&\le\mathcal E_\Omega((\mathcal Mv)\circ\Phi_h,
(\mathcal Mv)\circ\Phi_h)\\
&\le(1+C|h|)\mathcal E_{\Omega_h}(\mathcal Mv,\mathcal Mv)
\le(1+C|h|)\mathcal E_N(v,v).
\end{aligned}
\end{equation}
Then, by \eqref{e3.12} and \eqref{e3.13}, we conclude that
\begin{equation*}
\mathfrak a(T_hv,T_hv)=\mathcal{E}_N(T_hv,T_hv)+\|T_hv\|_{L^2(\Omega)}^2\leq(1+C|h|)\mathfrak a(v,v).
\end{equation*}
This finishes the proof of Lemma~\ref{l5}.
\end{proof}

Suppose now that $f\in(H^s(\Omega))^*$ and $v\in H^s(\Omega)$
satisfy $\mathfrak a(v,w)=f(w)$ for any $w\in H^s(\Omega)$.
By Lemma~\ref{l5}, we conclude that $d_h:=T_hv-v$ belongs to
$H^s(\Omega)$.
From the symmetry of $\mathfrak a$, we deduce that
\begin{align*}
\mathfrak a(d_h,d_h)
&=\mathfrak a(v+d_h,v+d_h)-\mathfrak a(v,v)-2\mathfrak a(v,d_h)\\
&=\mathfrak a(T_hv,T_hv)-\mathfrak a(v,v)-2f(d_h).
\end{align*}
By \eqref{e2.23} and Lemma~\ref{l5}, we find that
$\mathfrak a(d_h,d_h) \leq C|h|\mathfrak a(v,v)+2|f(d_h)|$.
From this and \eqref{e2.23}, it follows that
\begin{equation}\label{e3.14}
 \|d_h\|_{H^s(\Omega)}^2
 \leq C|h|\|v\|_{H^s(\Omega)}^2+C|f(d_h)|.
\end{equation}

\subsection{Difference estimates}\label{s3.2}

Let $\{\chi_j\}_{j=1}^{N_{\rm loc}}\subset C_{\rm c}^\infty(\mathbb{R}^n)$
be nonnegative functions forming a finite smooth partition of unity on a
neighborhood of $\overline\Omega$, subordinate to smaller boundary and
interior charts, and satisfying $\sum_{j=1}^{N_{\rm loc}}\chi_j=1$ there;
here $N_{\rm loc}$ denotes the total number of localization charts. For
each $j\in\{1,\ldots,N_{\rm loc}\}$, choose $\eta_j$ as in \eqref{e3.2} so
that $\eta_j\equiv1$ on a neighborhood of
$\operatorname{supp}\chi_j$. For any $h\in\mathcal C_j$, define
\begin{equation*}
\Phi_{j,h}(x):=x+\eta_j(x)h
\quad\text{and}\quad
T_{j,h}v:=v\circ\Phi_{j,h},
\end{equation*}
where $\mathcal C_j$ is the fixed open cone of admissible inward directions
for the corresponding chart. For interior charts, fix an arbitrary open cone
of admissible directions. Choose a common $h_0\in(0,1]$ for all these charts
such that the agreements in Section~\ref{s3.1} are admissible whenever
$0<|h|<h_0$.

Using the cone directions and repeated differences (see, for instance,
\cite[Propositions 2.1 and 2.2]{BN23} and
\cite[Propositions 2.1 and 2.2]{BLN23}), we obtain the following
characterization of Besov norms.

\begin{lemma}\label{l6}
Let $n\ge2$, let $s\in(0,1)$, and let $\Omega\subset\mathbb{R}^n$ be a bounded
Lipschitz domain. Let $r\in(0,1)$ and $v\in H^s(\Omega)$, and define
\begin{equation}\label{e3.17}
M_{s,r}(v):=\sum_{j=1}^{N_{\rm loc}}
\sup_{\substack{h\in\mathcal C_j\\0<|h|<h_0}}
|h|^{-r}\left\|T_{j,h}v-v\right\|_{H^s(\Omega)}.
\end{equation}
Then
\begin{equation*}
\|v\|_{B^{s+r}_{2,\infty}(\Omega)}
\le C\left[\|v\|_{H^s(\Omega)}+M_{s,r}(v)\right],
\end{equation*}
where the positive constant $C$ depends only on $n$, $s$, $r$, and $\Omega$.
\end{lemma}

\begin{proof}
From the openness of each admissible cone, we deduce that, for each
$j\in\{1,\ldots,N_{\rm loc}\}$, there exists a function $\theta_j$
satisfying
\begin{equation*}
\theta_j\in C_{\rm c}^\infty(\mathcal C_j\cap\{1/2<|z|<1\}),\quad
\theta_j\ge0,\quad\text{and}\quad
\int_{\mathbb{R}^n}\theta_j\,dz=1.
\end{equation*}
Choose the extension $Ev\in H^s(\mathbb{R}^n)$ from
Lemma~\ref{l1}\textup{(i)}. For any $\varepsilon\in(0,h_0)$, the operator
$A_\varepsilon$ is defined by setting, for any $x\in\mathbb{R}^n$,
\begin{equation*}
A_\varepsilon v(x):=\sum_{j=1}^{N_{\rm loc}}\chi_j(x)
\int_{\mathbb{R}^n}\theta_j(z)Ev(x+\varepsilon z)\,dz.
\end{equation*}
Recall that, for any $x\in\Omega\cap\operatorname{supp}\chi_j$ and
$z\in\operatorname{supp}\theta_j$, one has
$Ev(x+\varepsilon z)=T_{j,\varepsilon z}v(x)$. From the facts that
$\sum_j\chi_j=1$ and $\int_{\mathbb{R}^n}\theta_j\,dz=1$, it follows
that
\begin{equation*}
v(x)-A_\varepsilon v(x)=\sum_{j=1}^{N_{\rm loc}}\chi_j\int_{\mathbb{R}^n}\theta_j(z)
\left[v(x)-T_{j,\varepsilon z}v(x)\right]\,dz.
\end{equation*}
By Minkowski's inequality, we conclude that
\begin{equation}\label{e3.18}
\left\|v-A_\varepsilon v\right\|_{H^s(\Omega)}
\le C\sum_{j=1}^{N_{\rm loc}}\int_{\mathbb{R}^n}\theta_j(z)
\left\|v-T_{j,\varepsilon z}v\right\|_{H^s(\Omega)}\,dz
\le C\varepsilon^rM_{s,r}(v).
\end{equation}
From the integration by parts in $z$ and the fact that
$\int_{\mathbb{R}^n}\partial_i\theta_j\,dz=0$, we deduce that, for any
$i\in\{1,\ldots,n\}$,
\begin{align}\label{e3.19}
&\partial_i\left(\chi_j\int_{\mathbb{R}^n}
\theta_j(z)Ev(\,\cdot+\varepsilon z)\,dz\right)\nonumber\\
&\quad={}(\partial_i\chi_j)\int_{\mathbb{R}^n}\theta_j(z)
Ev(\,\cdot+\varepsilon z)\,dz
-\varepsilon^{-1}\chi_j\int_{\mathbb{R}^n}\partial_i\theta_j(z)
\left[T_{j,\varepsilon z}v-v\right]\,dz
\end{align}
on $\Omega$. By \eqref{e3.17}, we find that
\begin{subequations}
\begin{equation}\label{e3.20a}
\begin{aligned}
\left\|(\partial_i\chi_j)\int_{\mathbb{R}^n}\theta_j(z)Ev(\cdot+\varepsilon z)\,dz
-(\partial_i\chi_j)v\right\|_{H^s(\Omega)}\leq C\varepsilon^rM_{s,r}(v)
\end{aligned}
\end{equation}
and
\begin{align}\label{e3.20b}
\left\|\varepsilon^{-1}\chi_j\int_{\mathbb{R}^n}\partial_i\theta_j(z)
\left[T_{j,\varepsilon z}v-v\right]\,dz\right\|_{H^s(\Omega)}
&\leq C\varepsilon^{-1}\int|\partial_i\theta_j(z)|
|\varepsilon z|^rM_{s,r}(v)\,dz\notag\\
&\leq C\varepsilon^{r-1}M_{s,r}(v).
\end{align}
\end{subequations}
From \eqref{e2.3}, \eqref{e3.20a}, and \eqref{e3.20b}, we infer that
\begin{equation}\label{e3.21}
\|A_\varepsilon v\|_{H^{s+1}(\Omega)}
\le C\left[\|v\|_{H^s(\Omega)}
+\varepsilon^{r-1}M_{s,r}(v)\right].
\end{equation}
By the decomposition $v=(v-A_\varepsilon v)+A_\varepsilon v$,
\eqref{e3.18}, \eqref{e3.19}, \eqref{e3.20a}, \eqref{e3.20b}, and
\eqref{e3.21}, we conclude that, for any $\varepsilon\in(0,h_0)$,
\begin{equation}\label{e3.22}
\mathscr K(\varepsilon,v,H^s(\Omega),H^{s+1}(\Omega))
\le C\left[\varepsilon^rM_{s,r}(v)
+\varepsilon\|v\|_{H^s(\Omega)}\right].
\end{equation}
Using the decomposition $v=v+0$ for $\varepsilon\ge h_0$ and
\eqref{e3.22} for smaller $\varepsilon$, we find that
\begin{align*}
\|v\|_{B^{s+r}_{2,\infty}(\Omega)}
&\le C\sup_{\varepsilon>0}\varepsilon^{-r}
\mathscr K(\varepsilon,v,H^s(\Omega),H^{s+1}(\Omega))\\
&\le C\left[M_{s,r}(v)+\|v\|_{H^s(\Omega)}
+h_0^{-r}\|v\|_{H^s(\Omega)}\right]
\le C\left[M_{s,r}(v)+\|v\|_{H^s(\Omega)}\right].
\end{align*}
This finishes the proof of Lemma~\ref{l6}.
\end{proof}

\begin{lemma}\label{l7}
Let $n\ge2$, and let $\Omega\subset\mathbb{R}^n$ be a bounded Lipschitz
domain. Let $t,\sigma\in(0,1)$ and $v\in B^{t+\sigma}_{2,\infty}(\Omega)$.
Then, for any $j\in\{1,\ldots,N_{\rm loc}\}$ and any $h\in\mathcal C_j$
with $0<|h|<h_0$,
\begin{equation*}
\left\|T_{j,h}v-v\right\|_{H^t(\Omega)}\le C_{t,\sigma}|h|^\sigma
\|v\|_{B^{t+\sigma}_{2,\infty}(\Omega)},
\end{equation*}
where the positive constant $C_{t,\sigma}$ depends only on $n$, $t$,
$\sigma$, and the localization family, but is independent of $v$, $j$,
and $h$.
\end{lemma}

\begin{proof}
Extend $v$ by Lemma~\ref{l1}\textup{(i)} and let $V:=Ev$ and
$V_k:=\Delta_kV$ for any $k\in\mathbb{N}_0$, so that $V=\sum_{k\in\mathbb{N}_0}V_k$ is the
inhomogeneous dyadic decomposition of $V$ as in Section~\ref{s2.1}.
For any given $j\in\{1,\ldots,N_{\rm loc}\}$, let $\Phi_h:=\Phi_{j,h}$. From the change
of variables and \eqref{e3.4}, it follows that, for any $W\in H^t(\mathbb{R}^n)$ and any admissible $h$,
\begin{subequations}
\begin{equation}\label{e3.23a}
\|W\circ\Phi_h\|_{L^2(\mathbb{R}^n)}
\le C\|W\|_{L^2(\mathbb{R}^n)},
\end{equation}
\begin{equation}\label{e3.23b}
[W\circ\Phi_h]_{t,\mathbb{R}^n}^2
\le C_t[W]_{t,\mathbb{R}^n}^2,
\end{equation}
and
\begin{equation}\label{e3.23c}
\|W\circ\Phi_h\|_{H^t(\mathbb{R}^n)}
\le C_t\|W\|_{H^t(\mathbb{R}^n)}.
\end{equation}
\end{subequations}
By the fundamental theorem of calculus, we conclude that
\begin{equation}\label{e3.24}
V_k(\Phi_h(x))-V_k(x)
=\int_0^1\eta_j(x)h\cdot\nabla V_k(\Phi_{\tau h}(x))\,d\tau.
\end{equation}
From \eqref{e3.23a}--\eqref{e3.23c}, \eqref{e3.24}, Minkowski's
inequality, and the frequency support of $V_k$, we deduce that, for any
$k\in\mathbb{N}_0$,
\begin{align}\label{e3.25}
\|V_k\circ\Phi_h-V_k\|_{H^t(\mathbb{R}^n)}
&\le C\min\big\{\|V_k\|_{H^t(\mathbb{R}^n)},
|h|\|V_k\|_{H^{t+1}(\mathbb{R}^n)}\big\}\nonumber\\
&\le C\min\{1,2^k|h|\}\,2^{kt}\|V_k\|_{L^2(\mathbb{R}^n)}.
\end{align}
Choose $k_h\in\mathbb{N}_0$ such that
$2^{k_h}\le|h|^{-1}<2^{k_h+1}$. By the definition \eqref{equ2.3} of the
Besov norm, we have
$\|V_k\|_{L^2(\mathbb{R}^n)}
\le\|V\|_{B^{t+\sigma}_{2,\infty}(\mathbb{R}^n)}2^{-k(t+\sigma)}$.
Splitting the sum at $k_h$, we find that
\begin{equation}\label{e3.26}
\sum_{k\in\mathbb{N}_0}\min\{1,2^k|h|\}2^{-k\sigma}
\le|h|\sum_{k=0}^{k_h}2^{k(1-\sigma)}
+\sum_{k=k_h+1}^\infty2^{-k\sigma}
\le C_\sigma|h|^\sigma.
\end{equation}

Then, summing \eqref{e3.25} in $H^t(\mathbb{R}^n)$ and using \eqref{e3.26}, we
conclude that
\begin{equation*}
\|T_{j,h}v-v\|_{H^t(\Omega)}
\leq\|V\circ\Phi_h-V\|_{H^t(\mathbb{R}^n)}\leq C\|V\|_{B^{t+\sigma}_{2,\infty}(\mathbb{R}^n)}
\sum_{k\geq0}\min\{1,2^k|h|\}2^{-k\sigma}
\leq C|h|^\sigma\|v\|_{B^{t+\sigma}_{2,\infty}(\Omega)}.
\end{equation*}
This finishes the proof of Lemma~\ref{l7}.
\end{proof}

\subsection{Besov regularity and interpolation}\label{s3.3}

For any $t\ge0$, let $Y_t:=(H^t(\Omega))^*$ denote the dual space of
$H^t(\Omega)$, equipped with the norm
\begin{equation*}
\|f\|_{Y_t}:=\sup_{{\genfrac{}{}{0pt}{}{w\in H^t(\Omega)}
{\|w\|_{H^t(\Omega)}=1}}}|f(w)|.
\end{equation*}
From \eqref{e2.23}, we deduce that, for any $v,w\in H^s(\Omega)$,
\begin{equation*}
|\mathfrak a(v,w)|\le C\|v\|_{H^s(\Omega)}\|w\|_{H^s(\Omega)}
\quad\text{and}\quad
\mathfrak a(v,v)\ge c\|v\|_{H^s(\Omega)}^2.
\end{equation*}
For any $f\in Y_s$, the function $Sf$ is defined by, for any
$w\in H^s(\Omega)$,
\begin{equation}\label{e3.27}
\mathfrak a(Sf,w)=f(w).
\end{equation}
By the Lax--Milgram theorem and \eqref{e2.23}, we find that, for any
$f\in Y_s$, there exists a unique solution $Sf$, and the resulting inverse
operator $S$ is linear. From the fact that
$c\|Sf\|_{H^s(\Omega)}^2\le\mathfrak a(Sf,Sf)$, it follows that
\begin{equation*}
\|Sf\|_{H^s(\Omega)}\le C\|f\|_{Y_s}.
\end{equation*}

Following the regularity-improvement argument in
\cite[proof of Theorem 4.1]{BN23} and \cite[proof of Theorem 3.1]{BLN23},
we obtain the following assertion concerning the operator $S$.

\begin{proposition}\label{t3}
Let $n\ge2$, let $s\in(0,1)$, and let $\Omega\subset\mathbb{R}^n$ be a bounded
Lipschitz domain. Let $S$ be the inverse of the coercive form
$\mathfrak a$, as defined in \eqref{e3.27}, and $R\in(0,\min\{s,1/2\})$.
Then, for any $f\in(H^{s-R}(\Omega))^*$,
\begin{equation*}
\|Sf\|_{B^{s+R}_{2,\infty}(\Omega)}
\le C_R\|f\|_{(H^{s-R}(\Omega))^*}.
\end{equation*}
where the positive constant $C_R$ depends only on $n$, $s$, $R$, and
$\Omega$.
\end{proposition}

\begin{proof}
Let $v:=Sf$, $M_f:=\|f\|_{Y_{s-R}}$, and $t:=s-R>0$.
By the fact that $H^s(\Omega)\subset H^{s-R}(\Omega)$ and the energy
estimate, we deduce that
\begin{equation}\label{e3.28a}
\|v\|_{H^s(\Omega)}\le C\|f\|_{Y_s}\le CM_f.
\end{equation}
From the extension in Lemma~\ref{l1} and the definition \eqref{equ2.3} of
the Besov norm, it follows that
\begin{equation}\label{e3.28b}
\|v\|_{B^s_{2,\infty}(\Omega)}
\le C\|Ev\|_{B^s_{2,\infty}(\mathbb{R}^n)}
\le C\|Ev\|_{H^s(\mathbb{R}^n)}
\le C\|v\|_{H^s(\Omega)}\le CM_f.
\end{equation}

For any $k\in\mathbb{N}_0$, define the increasing regularity orders $\beta_k$ by
\begin{equation*}
\beta_0=0,\qquad \beta_{k+1}=\frac{R+\beta_k}{2},
\end{equation*}
and let $U_k:=\|v\|_{B^{s+\beta_k}_{2,\infty}(\Omega)}$. Then
$\beta_k=R(1-2^{-k})\to R$ as $k\to\infty$. Suppose that
$U_k<\infty$. By Lemma~\ref{l7} with $t:=s-R$ and
$\sigma:=R+\beta_k\in[R,2R]\subset(0,1)$, so that $t+\sigma=s+\beta_k$,
we conclude that, for any $j$ and any admissible $h$,
\begin{equation*}
\left|f(T_{j,h}v-v)\right|\le M_f\left\|T_{j,h}v-v\right\|_{H^{s-R}(\Omega)}
\le CM_f|h|^{R+\beta_k}U_k.
\end{equation*}
From \eqref{e3.14}, $2\beta_{k+1}=R+\beta_k<2R<1$, and $|h|\le1$, we
deduce that
\begin{equation*}
|h|^{-2\beta_{k+1}}\left\|T_{j,h}v-v\right\|_{H^s(\Omega)}^2
\le C|h|^{1-2\beta_{k+1}}M_f^2+CM_fU_k
\le C\left(M_f^2+M_fU_k\right).
\end{equation*}
By Lemma~\ref{l6}, \eqref{e3.28a}, and \eqref{e3.28b}, we find that
\begin{subequations}
\begin{equation}\label{e3.29a}
M_{s,\beta_{k+1}}(v)\le C\left(M_f^2+M_fU_k\right)^{1/2}
\end{equation}
and
\begin{equation}\label{e3.29b}
U_{k+1}^2
\le C\left[\|v\|_{H^s(\Omega)}^2+M_{s,\beta_{k+1}}(v)^2\right]
\le C\left(M_f^2+M_fU_k\right),
\end{equation}
where $M_{s,\beta_{k+1}}(v)$ is as in \eqref{e3.17}.
\end{subequations}

Since $\sigma\in[R,2R]$, $\beta_{k+1}\in[R/2,R]$, and $2R<1$, the
constants appearing in \eqref{e3.29a} and \eqref{e3.29b} are independent of
$k$. Choose $A$ such that $U_0\le AM_f$ and $C(1+A)\le A^2$. By
\eqref{e3.29a}, \eqref{e3.29b}, and induction, we conclude that, for any
$k\in\mathbb{N}_0$, if $U_k\le AM_f$, then
\begin{equation*}
U_{k+1}^2\le C(1+A)M_f^2\le A^2M_f^2,
\end{equation*}
and hence $\sup_{k\in\mathbb{N}_0}U_k\le AM_f$. From this, \eqref{e3.29a}, and
\eqref{e3.22}, we deduce that, for any fixed $\varepsilon\in(0,\infty)$ and any
$k\in\mathbb{N}$,
\begin{equation}\label{e3.30}
\mathscr K(\varepsilon,v,H^s(\Omega),H^{s+1}(\Omega))
\le CM_f\varepsilon^{\beta_k}.
\end{equation}
Letting $k\to\infty$ in \eqref{e3.30} at fixed $\varepsilon\in(0,\infty)$ and using
\eqref{e2.7}, we find that
\begin{equation*}
\mathscr K(\varepsilon,v,H^s(\Omega),H^{s+1}(\Omega))
\le CM_f\lim_{k\to\infty}\varepsilon^{\beta_k}
=CM_f\varepsilon^R,
\end{equation*}
and hence
\begin{equation*}
\|v\|_{B^{s+R}_{2,\infty}(\Omega)}
\le C\sup_{\varepsilon>0}\varepsilon^{-R}
\mathscr K(\varepsilon,v,H^s(\Omega),H^{s+1}(\Omega))
\le CM_f.
\end{equation*}
This finishes the proof of Proposition~\ref{t3}.
\end{proof}

\begin{proposition}\label{t4}
Let $n\geq2$, let $s\in(0,1)$, and let $\Omega\subset\mathbb{R}^n$ be a bounded
Lipschitz domain. Let $S$ be the inverse of the coercive form
$\mathfrak a$, as defined in \eqref{e3.27}.
Then, for any $f\in(H^{s/2}(\Omega))^*$,
\begin{equation*}
\|Sf\|_{H^{3s/2}(\Omega)}\leq C\|f\|_{(H^{s/2}(\Omega))^*},
\end{equation*}
where the positive constant $C$ depends only on $n$, $s$, and $\Omega$.
\end{proposition}

\begin{proof}
Choose $R$ to be the midpoint of the interval
$(s/2,\min\{s,1/2\})$ and let $\theta:=s/(2R)\in(0,1)$. From the energy
estimate and Proposition~\ref{t3}, it follows that both
$S:Y_s\to H^s(\Omega)$ and
$S:Y_{s-R}\to B^{s+R}_{2,\infty}(\Omega)$ are bounded. By the embedding
$Y_{s-R}\subset Y_s$ and uniqueness in \eqref{e3.27}, we conclude that
these maps agree on their common domain.

From the extension and restriction operators in
Lemma~\ref{l1}\textup{(i)}, we deduce that, for any $a\in[s-R,s]$,
\begin{equation}\label{e3.31a}
\mathcal R^*:Y_a\longrightarrow H^{-a}(\mathbb{R}^n)
\quad\text{and}\quad
E^*:H^{-a}(\mathbb{R}^n)\longrightarrow Y_a.
\end{equation}
By this, we find that, for any $f\in Y_a$, $V\in H^a(\mathbb{R}^n)$,
$G\in H^{-a}(\mathbb{R}^n)$, and $w\in H^a(\Omega)$,
\begin{equation}\label{e3.31b}
(\mathcal R^*f)(V)=f(\mathcal RV)
\quad\text{and}\quad
(E^*G)(w)=\langle G,Ew\rangle.
\end{equation}
From $\mathcal RE=I$, it follows that
\begin{equation}\label{e3.33}
E^*\mathcal R^*=(\mathcal RE)^*=I.
\end{equation}
Let $a_\theta:=s-\theta R$. By the full-space interpolation, we conclude that
\begin{equation*}
(H^{-s}(\mathbb{R}^n),H^{-s+R}(\mathbb{R}^n))_{\theta,2}
=H^{-a_\theta}(\mathbb{R}^n)
\end{equation*}
with equivalent norms. From this and the interpolation of \eqref{e3.31a},
 together with \eqref{e3.31b} and \eqref{e3.33}, we infer that
\begin{equation*}
\|f\|_{Y_{a_\theta}}\le C\|\mathcal R^*f\|_{H^{-a_\theta}(\mathbb{R}^n)}
\le C\|f\|_{(Y_s,Y_{s-R})_{\theta,2}}
\end{equation*}
and
\begin{equation*}
\|f\|_{(Y_s,Y_{s-R})_{\theta,2}}
\le C\|\mathcal R^*f\|_{H^{-a_\theta}(\mathbb{R}^n)}
\le C\|f\|_{Y_{a_\theta}}.
\end{equation*}
By this, we find that
\begin{equation}\label{equ3.30}
(Y_s,Y_{s-R})_{\theta,2}=Y_{s-\theta R}=Y_{s/2}
\end{equation}
with equivalent norms. Moreover, from Lemma~\ref{l1}\textup{(iv)}, it follows that
\begin{equation*}
\left(H^s(\Omega),B^{s+R}_{2,\infty}(\Omega)\right)_{\theta,2}
=H^{s+\theta R}(\Omega)=H^{3s/2}(\Omega)
\end{equation*}
with equivalent norms. By this and \eqref{equ3.30}, we conclude that, for
any $f\in Y_{s/2}$,
\begin{equation*}
\|Sf\|_{H^{3s/2}(\Omega)}\le C\|f\|_{Y_{s/2}}.
\end{equation*}
This finishes the proof of Proposition~\ref{t4}.
\end{proof}

Let $q\in L^2(\Omega^c,\sigma_s)$ satisfy $\int_{\Omega^c}q\,d\sigma_s=0$ and
$u$ be the weak solution to problem \eqref{e1.1} with $g=a_sq$.
For $w\in H^{s/2}(\Omega)$, define
\begin{equation*}
 f(w):=F_q(w)+(u,w)_{L^2(\Omega)}.
\end{equation*}
From Lemma~\ref{l4}\textup{(iii)}, we deduce that
$\mathfrak a(u,w)=f(w)$ for any $w\in H^s(\Omega)$.
By Lemma~\ref{l4}\textup{(iv)}, we find that
\begin{equation*}
|f(w)|\leq\big(C\|q\|_{L^2(\Omega^c,\sigma_s)}+\|u\|_{L^2(\Omega)}\big)
\|w\|_{H^{s/2}(\Omega)}.
\end{equation*}
From \eqref{e2.19}, it follows that
$\|f\|_{Y_{s/2}}\leq C\|q\|_{L^2(\Omega^c,\sigma_s)}$.
By Proposition~\ref{t4}, we conclude that
\begin{equation}\label{e3.34}
\|u\|_{H^{3s/2}(\Omega)}\leq C\|q\|_{L^2(\Omega^c,\sigma_s)}.
\end{equation}

\section{Exterior reflection difference}\label{s4}

Throughout this section, let $n\geq2$, $s\in(0,1)$, and $a:=3s/2$.
In this section, by means of a reflection argument, we establish the exterior reflection difference bound for
the weak solution to \eqref{e1.1}.

\subsection{Reflected difference estimate below first order}\label{s4.1}

For any $X:=(x',t),Y:=(y',r)\in\mathbb{R}^n_+$, let
\begin{equation*}
K_{\rm ref}(X,Y):=\big(|x'-y'|^2+(t+r)^2\big)^{-(n+2s)/2}.
\end{equation*}
For any measurable function $w:\mathbb{R}^n_+\to\mathbb{R}$, the
\emph{reflected difference seminorm} $J_s(w)$ is defined by setting
\begin{equation}\label{e4.1}
J_s(w):=\left[\iint_{\mathbb{R}^n_+\times\mathbb{R}^n_+}
r^{-s}K_{\rm ref}(X,Y)|w(X)-w(Y)|^2\,dX\,dY\right]^{1/2},
\end{equation}
where $dX:=dx'\,dt$ and $dY:=dy'\,dr$.

\begin{lemma}\label{l8}
Let $n\ge2$ and $s\in(0,1)$. Fix $\alpha\in(0,\min\{1,s+1/2\})$ and let
$\lambda:=2s-2\alpha$. Then, for any $w\in H^\alpha(\mathbb{R}^n_+)$,
\begin{equation}\label{equ4.1}
\iint_{\mathbb{R}^n_+\times\mathbb{R}^n_+}
r^\lambda K_{\rm ref}(X,Y)|w(X)-w(Y)|^2\,dX\,dY
\le C_{n,s,\alpha}\|w\|_{H^\alpha(\mathbb{R}^n_+)}^2,
\end{equation}
where the positive constant $C_{n,s,\alpha}$ depends only on $n$, $s$,
and $\alpha$.
\end{lemma}

\begin{remark}
  When $\lambda=0$, \eqref{equ4.1} follows directly from even extension. Lemma \ref{l8} is its weighted version.
\end{remark}

\begin{proof}[Proof of Lemma \ref{l8}]
Denote by $\mathcal C_{s,\alpha}(w)$ the square root of the left-hand side
of \eqref{equ4.1}. Since $\alpha<s+1/2$, we find that $\lambda>-1$.

We first consider the case $0<\alpha<1/2$. For any $t,r\in(0,\infty)$,
\begin{equation}\label{e4.2}
\int_0^\infty r^\lambda(t+r)^{-1-2s}\,dr
=\mathrm B(\lambda+1,2\alpha)\,t^{-2\alpha}
\end{equation}
and
\begin{equation}\label{e4.3}
r^\lambda\int_0^\infty(t+r)^{-1-2s}\,dt
=\frac{1}{2s}\,r^{-2\alpha},
\end{equation}
where, for any $a,b\in(0,\infty)$, $\mathrm B(a,b)$ denotes the beta function, i.e.,
$\mathrm B(a,b):=\int_0^\infty q^{a-1}(1+q)^{-a-b}\,dq$.
By the change of variables, we find that, for any $t,r\in(0,\infty)$,
\begin{equation}\label{e4.4}
\int_{\mathbb{R}^{n-1}}\big[|z|^2+(t+r)^2\big]^{-(n+2s)/2}\,dz
=\frac{\pi^{(n-1)/2}\Gamma(s+1/2)}{\Gamma(n/2+s)}(t+r)^{-1-2s}.
\end{equation}
From \eqref{e4.2}, \eqref{e4.3}, \eqref{e4.4}, the fact that
$|w(X)-w(Y)|^2\le2|w(X)|^2+2|w(Y)|^2$, and Hardy's inequality in
Lemma~\ref{l1}(iii), we deduce that
\begin{equation*}
\mathcal C_{s,\alpha}(w)^2
\le C\int_{\mathbb{R}^n_+}t^{-2\alpha}|w(x',t)|^2\,dx'\,dt
\le C\|w\|_{H^\alpha(\mathbb{R}^n_+)}^2.
\end{equation*}

Next, we consider the case $1/2<\alpha<\min\{1,s+1/2\}$. Let
$h:=\operatorname{Tr}w\in H^{\alpha-1/2}(\mathbb{R}^{n-1})$ and
$\kappa:=(1+|\xi'|^2)^{1/2}$. For any $\xi'\in\mathbb{R}^{n-1}$ and
$t\in(0,\infty)$, define $Ph$ by setting
\begin{equation}\label{e4.5}
\mathcal F_{x'}(Ph)(\xi',t):=e^{-\kappa t}\mathcal F_{x'}h(\xi'),
\end{equation}
where $\mathcal F_{x'}$ denotes the Fourier transform with respect to $x'$.
Let $\widetilde{Ph}(x',t):=Ph(x',|t|)$. From the Fourier transform on
$\mathbb{R}^n$ and Plancherel's theorem, it follows that, for any
$\xi'\in\mathbb{R}^{n-1}$ and $\tau\in\mathbb{R}$,
\begin{equation}\label{e4.6a}
\widehat{\widetilde{Ph}}(\xi',\tau)
=\sqrt{2/\pi}\,\frac{\kappa}{\kappa^2+\tau^2}\mathcal F_{x'}h(\xi')
\end{equation}
and
\begin{equation}\label{e4.6b}
\|Ph\|_{H^\alpha(\mathbb{R}^n_+)}^2
\le\|\widetilde{Ph}\|_{H^\alpha(\mathbb{R}^n)}^2
\le C\int_{\mathbb{R}^{n-1}}\kappa^{2\alpha-1}
|\mathcal F_{x'}h(\xi')|^2\,d\xi'.
\end{equation}

Let $w_0:=w-Ph$. From \eqref{e2.5}, \eqref{e2.6}, \eqref{e4.2},
\eqref{e4.3}, \eqref{e4.4}, \eqref{e4.5}, \eqref{e4.6a}, and
\eqref{e4.6b}, we deduce that
\begin{equation}\label{e4.7a}
\operatorname{Tr}w_0=0,\quad
\|w_0\|_{H^\alpha(\mathbb{R}^n_+)}
\le C\|w\|_{H^\alpha(\mathbb{R}^n_+)},
\end{equation}
and
\begin{equation}\label{e4.7c}
\mathcal C_{s,\alpha}(w_0)^2
\le C\int_{\mathbb{R}^n_+}t^{-2\alpha}|w_0(x',t)|^2\,dx'\,dt
\le C\|w_0\|_{H^\alpha(\mathbb{R}^n_+)}^2.
\end{equation}

To estimate $\mathcal C_{s,\alpha}(Ph)$, let
\begin{equation*}
b_{n,s}:=\left[\int_{\mathbb{R}^{n-1}}
(1+|z|^2)^{-(n+2s)/2}\,dz\right]^{-1},
\quad
k_s(z):=b_{n,s}(1+|z|^2)^{-(n+2s)/2},
\end{equation*}
and, for any $\tau\in[0,\infty)$, define
\begin{equation*}
\Phi_s(\tau):=\int_{\mathbb{R}^{n-1}}k_s(z)\cos(\tau z_1)\,dz.
\end{equation*}
Since $0\le1-\cos q\le2|q|$ for any $q\in\mathbb{R}$, it follows that, for
any $\tau\in[0,\infty)$, $|\Phi_s(\tau)|\le1$ and
\begin{equation}\label{e4.8}
0\le1-\Phi_s(\tau)\le C_{n,s}\tau.
\end{equation}
For any $\vartheta\in[0,1]$, let
\begin{equation*}
I(\vartheta):=\int_0^\infty\int_0^\infty\frac{r^\lambda}{(t+r)^{1+2s}}
\left[e^{-2t}+e^{-2r}-2\Phi_s(\vartheta(t+r))e^{-(t+r)}\right]\,dt\,dr.
\end{equation*}
By Plancherel's theorem for $Ph$ on $\mathbb{R}^{n-1}$ and the change of
variables $(t,r)\mapsto(\kappa t,\kappa r)$, we conclude that
\begin{equation}\label{e4.9}
\mathcal C_{s,\alpha}(Ph)^2
=b_{n,s}^{-1}\int_{\mathbb{R}^{n-1}}\kappa^{2\alpha-1}
I(|\xi'|/\kappa)\left|\widehat h(\xi')\right|^2\,d\xi'.
\end{equation}
Let $R:=t+r$ and $z:=r/R\in[0,1]$, so that $r=Rz$ and $t=R(1-z)$. From this
change of variables, we deduce that
\begin{align}\label{e4.10a}
I(\vartheta)
&=\int_0^\infty\int_0^1 R^{-2\alpha}z^\lambda
\left\{(e^{-t}-e^{-r})^2+2[1-\Phi_s(\vartheta R)]e^{-R}\right\}\,dz\,dR
\nonumber\\
&=\int_0^1\int_0^1 R^{-2\alpha}z^\lambda
\left\{(e^{-t}-e^{-r})^2+2[1-\Phi_s(\vartheta R)]e^{-R}\right\}\,dz\,dR
\nonumber\\
&\quad+\int_1^\infty\int_0^1\cdots
=:\mathrm I_1+\mathrm I_2.
\end{align}
Moreover, by \eqref{e4.8} and the fact that
$|e^{-t}-e^{-r}|\le|t-r|\le R$, we find that, for any $R\in(0,1)$ and
$\vartheta\in[0,1]$,
\begin{equation*}
(e^{-t}-e^{-r})^2+2\left[1-\Phi_s(\vartheta R)\right]e^{-R}\le C_sR,
\end{equation*}
which, combined with \eqref{e4.10a}, implies that
\begin{equation}\label{e4.11}
\mathrm I_1
\le C\int_0^1R^{1-2\alpha}\,dR\int_0^1z^\lambda\,dz
=\frac{C}{(2-2\alpha)(\lambda+1)}<\infty.
\end{equation}
From the fact that $|\Phi_s|\le1$, it follows that, for any
$R\in[1,\infty)$,
\begin{equation*}
(e^{-t}-e^{-r})^2+2[1-\Phi_s(\vartheta R)]e^{-R}
\le2(e^{-2t}+e^{-2r}),
\end{equation*}
which further implies that, for any $R\in[1,\infty)$,
\begin{equation*}
\int_0^1z^\lambda e^{-2(1-z)R}\,dz
\le e^{-R}\int_0^{1/2}z^\lambda\,dz
+C_\lambda\int_{1/2}^1e^{-2(1-z)R}\,dz
\le C_\lambda R^{-1}
\end{equation*}
and
\begin{equation*}
\int_0^1z^\lambda e^{-2zR}\,dz
=R^{-\lambda-1}\int_0^R q^\lambda e^{-2q}\,dq
\le C_\lambda R^{-\lambda-1}.
\end{equation*}
By these estimates, we conclude that
\begin{equation}\label{e4.14}
{\rm I}_2\le C\int_1^\infty R^{-2\alpha-1}\,dR+C\int_1^\infty R^{-2s-1}\,dR<\infty.
\end{equation}
Then, from \eqref{e4.9}, \eqref{e4.10a}, \eqref{e4.11}, \eqref{e4.14}, and the trace inequality
\eqref{e2.6}, it follows that
\begin{equation*}
\mathcal C_{s,\alpha}(Ph)^2\le C\int_{\mathbb{R}^{n-1}}\kappa^{2\alpha-1}
\left|\widehat h(\xi')\right|^2\,d\xi'
\le C\|w\|_{H^\alpha(\mathbb{R}^n_+)}^2,
\end{equation*}
which, combined with \eqref{e4.7a} and \eqref{e4.7c}, implies that
\begin{equation*}
\mathcal C_{s,\alpha}(w)
\le\mathcal C_{s,\alpha}(w_0)+\mathcal C_{s,\alpha}(Ph)
\le C\|w\|_{H^\alpha(\mathbb{R}^n_+)}.
\end{equation*}

Finally, we consider the case $\alpha=1/2$. Fix
$\varepsilon\in(0,\min\{s,1/2\})$ and let $\alpha_\pm:=1/2\pm\varepsilon$. Let
$\mathcal Ew(x',t):=w(x',|t|)$. By this definition, we conclude that,
for any $\beta\in(0,1)$,
\begin{equation}\label{e4.16a}
\|\mathcal Ew\|_{L^2(\mathbb{R}^n)}^2
=2\|w\|_{L^2(\mathbb{R}^n_+)}^2
\quad\text{and}\quad
[\mathcal Ew]_{\beta,\mathbb{R}^n}^2
\le4[w]_{\beta,\mathbb{R}^n_+}^2.
\end{equation}
From the even extension estimates \eqref{e4.16a}, we deduce that
\begin{equation}\label{equ4.17}
\left(H^{\alpha_-}(\mathbb{R}^n_+),H^{\alpha_+}(\mathbb{R}^n_+)\right)_{1/2,2}
=H^{1/2}(\mathbb{R}^n_+)
\end{equation}
with equivalent norms. Let $W_\pm:=r^{2s-2\alpha_\pm}K_{\rm ref}$. By
the interpolation of weighted $L^2$ spaces, we find that
\begin{equation*}
\left(L^2(\mathbb{R}^n_+\times\mathbb{R}^n_+,W_-\,dX\,dY),
L^2(\mathbb{R}^n_+\times\mathbb{R}^n_+,W_+\,dX\,dY)\right)_{1/2,2}
=L^2(\mathbb{R}^n_+\times\mathbb{R}^n_+,\sqrt{W_-W_+}\,dX\,dY)
\end{equation*}
with equivalent norms. Indeed, define $\mathscr K_2$ by replacing the sum
of norms in $\mathscr K$ by the square root of the sum of their squares.
Taking the pointwise minimization, we find that, for any
$\tau\in(0,\infty)$,
\begin{equation*}
\mathscr K_2(\tau,F)^2
=\iint_{\mathbb{R}^n_+\times\mathbb{R}^n_+}
\frac{\tau^2W_-W_+}{W_-+\tau^2W_+}|F|^2\,dX\,dY.
\end{equation*}
By Fubini's theorem and $\alpha_-+\alpha_+=1$, we conclude that
\begin{align*}
\int_0^\infty\tau^{-2}\mathscr K_2(\tau,F)^2\,d\tau
&=\frac{\pi}{2}\iint_{\mathbb{R}^n_+\times\mathbb{R}^n_+}
\sqrt{W_-W_+}|F|^2\,dX\,dY\\
&=\frac{\pi}{2}\iint_{\mathbb{R}^n_+\times\mathbb{R}^n_+}
r^{2s-1}K_{\rm ref}(X,Y)|F|^2\,dX\,dY.
\end{align*}
From this, the interpolation \eqref{equ4.17}, and the results in the case
$\alpha\ne1/2$, we deduce that
\begin{equation*}
\iint_{\mathbb{R}^n_+\times\mathbb{R}^n_+}
r^{2s-1}K_{\rm ref}(X,Y)|w(X)-w(Y)|^2\,dX\,dY
\le C\|w\|_{H^{1/2}(\mathbb{R}^n_+)}^2.
\end{equation*}
This finishes the proof of Lemma~\ref{l8}.
\end{proof}

Applying Lemma~\ref{l8} with $\alpha:=3s/2$ and $\lambda:=-s$, which is
admissible whenever $s\in(0,2/3)$ because
$\alpha=3s/2<\min\{1,s+1/2\}$, we find that, for any
$w\in H^{3s/2}(\mathbb{R}^n_+)$,
\begin{equation}\label{e4.17}
J_s(w)^2\le C_{n,s}\|w\|_{H^{3s/2}(\mathbb{R}^n_+)}^2.
\end{equation}

\subsection{First-order weighted estimate}\label{s4.2}

\begin{lemma}\label{l9}
Let $n\ge2$, $s\in(0,1)$, and $w\in H^1(\mathbb{R}^n_+)$. Then the
reflected difference quantity $J_s(w)$ defined in \eqref{e4.1} satisfies
\begin{equation*}
J_s(w)^2\le C_{n,s}\int_{\mathbb{R}^n_+}t^{2-3s}
|\nabla w(x',t)|^2\,dx'\,dt,
\end{equation*}
where the positive constant $C_{n,s}$ depends only on $n$ and $s$.
\end{lemma}

\begin{proof}
From Fubini's theorem and the identity
$q^{-n-2s}=(n+2s)\int_q^\infty R^{-n-2s-1}\,dR$ for any $q\in(0,\infty)$,
it follows that
\begin{equation*}
J_s(w)^2\le C\int_0^\infty R^{-n-2s-1}\int_0^R\int_0^R r^{-s}\int_{\mathbb{R}^{n-1}}\int_{|z|<R}
|w(x',t)-w(x'+z,r)|^2\,dz\,dx'\,dr\,dt\,dR.
\end{equation*}
Fix $0<t,r<R$. Then
\begin{align}\label{equ4.19}
w(x',t)-w(x'+z,r)
&=[w(x',t)-w(x',R)]+[w(x',R)-w(x'+z,R)]\notag\\
&\quad+[w(x'+z,R)-w(x'+z,r)].
\end{align}
By the fundamental theorem of calculus, we find that, for any
$R\in(0,\infty)$ and any $z\in\mathbb{R}^{n-1}$,
\begin{equation*}
\int_{\mathbb{R}^{n-1}}|w(x',R)-w(x'+z,R)|^2\,dx'
\le|z|^2\|\nabla'w(\,\cdot,R)\|_{L^2(\mathbb{R}^{n-1})}^2,
\end{equation*}
where $\nabla'$ denotes the gradient with respect to $x'$. From the
integration in $z$, $t$, and $r$, we deduce that the horizontal
contribution is bounded by
\begin{align}\label{e4.18}
&C\int_0^\infty R^{-n-2s-1}R^{n+1}R\frac{R^{1-s}}{1-s}
\|\nabla'w(\,\cdot,R)\|_{L^2(\mathbb{R}^{n-1})}^2\,dR\notag\\
&\quad=C_{n,s}\int_0^\infty R^{2-3s}
\|\nabla'w(\,\cdot,R)\|_{L^2(\mathbb{R}^{n-1})}^2\,dR.
\end{align}

By the identity $\int_{|z|<R}dz=C_nR^{n-1}$, we conclude that the vertical contributions are bounded by
constant multiples of
\begin{equation*}
 {\rm I}_3:=\int_0^\infty R^{-1-3s}\int_0^R\|w(\cdot,t)-w(\cdot,R)\|_{L^2(\mathbb{R}^{n-1})}^2\,dt\,dR
\end{equation*}
and
\begin{equation*}
 {\rm I}_4:=\int_0^\infty R^{-1-2s}\int_0^Rr^{-s}\|w(\cdot,r)-w(\cdot,R)\|_{L^2(\mathbb{R}^{n-1})}^2\,dr\,dR
\end{equation*}
respectively.
From the identity $w(\cdot,R)-w(\cdot,t)=\int_t^R\partial_\tau w(\cdot,\tau)\,d\tau$, we deduce that, for any $0<t<R$,
\begin{align}\label{e4.19}
\|w(\cdot,R)-w(\cdot,t)\|_{L^2(\mathbb{R}^{n-1})}^2
&\leq\left(\int_t^R\tau^{s-1}\,d\tau\right)
\int_t^R\tau^{1-s}\|(\partial_tw)(\cdot,\tau)\|_{L^2(\mathbb{R}^{n-1})}^2\,d\tau\nonumber\\
&\leq\frac{R^s}{s}\int_t^R\tau^{1-s}\|(\partial_tw)(\cdot,\tau)\|_{L^2(\mathbb{R}^{n-1})}^2\,d\tau.
\end{align}
By Fubini's theorem applied to \eqref{e4.19}, we find that
\begin{align}\label{e4.20}
{\rm I}_3&\leq\frac1s\int_0^\infty R^{-1-2s}
\int_0^R\tau^{2-s}\|\partial_tw(\cdot,\tau)\|_{L^2(\mathbb{R}^{n-1})}^2\,d\tau\,dR\nonumber\\
&=\frac1s\int_0^\infty\tau^{2-s}\|\partial_tw(\cdot,\tau)\|_{L^2(\mathbb{R}^{n-1})}^2
\int_\tau^\infty R^{-1-2s}\,dR\,d\tau\nonumber\\
&\leq\frac1{2s^2}\int_0^\infty\tau^{2-3s}\|(\partial_tw)(\cdot,\tau)\|_{L^2(\mathbb{R}^{n-1})}^2\,d\tau.
\end{align}
From the identity $\int_0^\tau r^{-s}\,dr=\tau^{1-s}/(1-s)$, it follows that
\begin{align}\label{e4.21}
{\rm I}_4&\leq\frac1{s(1-s)}\int_0^\infty\tau^{2-2s}\|\partial_tw(\cdot,\tau)\|_{L^2(\mathbb{R}^{n-1})}^2
\int_\tau^\infty R^{-1-s}\,dR\,d\tau\nonumber\\
&\leq\frac1{s^2(1-s)}\int_0^\infty\tau^{2-3s}
\|(\partial_tw)(\cdot,\tau)\|_{L^2(\mathbb{R}^{n-1})}^2\,d\tau.
\end{align}
Finally, by \eqref{equ4.19}, \eqref{e4.18}, \eqref{e4.20}, and
\eqref{e4.21}, we conclude that
\begin{equation*}
J_s(w)^2\leq C\left(\int_0^\infty R^{2-3s}
\|\nabla'w(\cdot,R)\|_{L^2(\mathbb{R}^{n-1})}^2\,dR+{\rm I}_3+{\rm I}_4\right)
\leq C\int_{\mathbb{R}^n_+}t^{2-3s}|\nabla w|^2\,dx'\,dt.
\end{equation*}
This finishes the proof of Lemma~\ref{l9}.
\end{proof}

\subsection{Reflection difference}\label{s4.3}

Let $\Omega\subset\mathbb{R}^n$ be a bounded Lipschitz domain. Choose a
finite boundary atlas with $N_0$ charts, indexed by
$j\in\{1,\ldots,N_0\}$, and fix four nested neighborhoods in each chart,
\begin{equation*}
U_j^0\Subset U_j^1\Subset U_j^2\Subset U_j^3,
\end{equation*}
such that the sets $U_j^0$ still cover $\partial\Omega$. Shrinking the
inner neighborhoods if necessary, we may assume that the vertical
reflections of $U_j^1\cap\Omega^c$ lie in $U_j^2\cap\Omega$, and that the
largest neighborhoods are contained in the original charts. The
neighborhoods are chosen small enough that their active exterior portions
satisfy $\delta(z)<r_*$, with $r_*$ as in Lemma~\ref{l2}.

In the coordinates of chart $j$, let $\varphi_j$ denote its graph
function, extended to $\mathbb{R}^{n-1}$ with Lipschitz constant at most
$L$. For any $x'\in\mathbb{R}^{n-1}$ and $t\in[0,\infty)$, define
\begin{equation*}
F_j(x',t):=(x',t+\varphi_j(x')),
\end{equation*}
and, for any $y'\in\mathbb{R}^{n-1}$ and $r\in(0,\infty)$, let
\begin{equation*}
R_j\big(F_j(y',-r)\big):=F_j(y',r).
\end{equation*}
The map $F_j$ takes flat coordinates to graph coordinates, and $R_j$
reflects an exterior point vertically across the graph. From the Lipschitz
bound for $\varphi_j$, it follows that, for any
$x',y'\in\mathbb{R}^{n-1}$ and $t,r\in(0,\infty)$,
\begin{equation}\label{e4.25}
|F_j(x',t)-F_j(y',-r)|\sim\left[|x'-y'|^2+(t+r)^2\right]^{1/2},
\end{equation}
and that, for any $x=F_j(x',t)$ in the chart,
\begin{equation}\label{e4.26}
\delta(F_j(x',t))\sim|t|.
\end{equation}
For any $x',y'\in\mathbb{R}^{n-1}$ and $t,r\in(0,\infty)$, let
$d:=|F_j(x',t)-F_j(y',-r)|$. By the triangle inequality, we conclude that
$|x'-y'|\le d$,
\begin{equation*}
t+r\le|t+r+\varphi_j(x')-\varphi_j(y')|+L|x'-y'|\le(1+L)d,
\end{equation*}
and
\begin{equation*}
d\le(1+L)|x'-y'|+t+r.
\end{equation*}
Choose $\psi_j\in C_{\rm c}^\infty(U_j^3)$ such that
$0\le\psi_j\le1$, $\psi_j=1$ on $U_j^2$, and
$\|\partial^\alpha\psi_j\|_{L^\infty(\mathbb{R}^n)}
\le C_\alpha r_0^{-|\alpha|}$ for any $|\alpha|\le2$,
with $C_\alpha$ controlled by the fixed Lipschitz character. For any
$v\in H^{3s/2}(\Omega)$, take an ordinary extension
$Ev\in H^a(\mathbb{R}^n)$ and define, for any $X\in\mathbb{R}^n_+$,
\begin{equation}\label{e4.27}
w_j(X):=(\psi_jEv)(F_j(X)).
\end{equation}
From the fact that
$F_j(\mathbb{R}^n_+)\cap\operatorname{supp}\psi_j\subset\Omega$, it follows
that $w_j$ is independent of the exterior values of $Ev$.

\begin{lemma}\label{l10}
Let $n\ge2$, let $s\in(0,1)$, and let $\Omega\subset\mathbb{R}^n$ be a bounded
Lipschitz domain. For any $v\in H^{3s/2}(\Omega)$, define $w_j$ by
\eqref{e4.27} for each $j\in\{1,\ldots,N_0\}$. Then
\begin{equation}\label{e4.29}
\sum_{j=1}^{N_0}J_s(w_j)^2
\le C\|v\|_{H^{3s/2}(\Omega)}^2,
\end{equation}
where the positive constant $C$ depends only on $n$, $s$, and $\Omega$.
\end{lemma}

\begin{proof}
By the bi-Lipschitz change of variables in \eqref{e4.27}, we find that
$\|w_j\|_{H^{3s/2}(\mathbb{R}^n_+)}\le C\|v\|_{H^{3s/2}(\Omega)}$
whenever $s\in(0,2/3)$. From \eqref{e4.17} and the finite number of
charts, it follows that \eqref{e4.29} holds in this case.

For $s\in[2/3,1)$, let $\gamma:=3s/2-1\in[0,1/2)$. By
Lemma~\ref{l1}(ii) and \eqref{e4.17}, we conclude that
$w_j\in H^1(\mathbb{R}^n_+)$ and that $w_j$ has bounded support. Let
$f:=\psi_jv$. Then, for any $i\in\{1,\ldots,n-1\}$ and any
$X\in\mathbb{R}^n_+$,
\begin{equation*}
\partial_iw_j
=(\partial_if)\circ F_j+(\partial_i\varphi_j)(\partial_nf)\circ F_j
\quad\text{and}\quad
\partial_tw_j=(\partial_nf)\circ F_j,
\end{equation*}
which, combined with the facts that
$\|DF_j\|_{L^\infty(\mathbb{R}^n_+)}\le C(1+L)$ and $t\sim\delta(F_j(x',t))$ on the support of $w_j$, implies that
\begin{equation*}
\int_{\mathbb{R}^n_+}t^{-2\gamma}|\nabla w_j(x',t)|^2\,dx'\,dt
\le C\int_{U_j^3\cap\Omega}\delta^{-2\gamma}
\big(|\nabla v|^2+|v|^2\big)\,dx.
\end{equation*}
From this, Hardy's inequality for $v$ and $\nabla v$, and
Lemma~\ref{l9}, we deduce that
\begin{equation*}
\sum_{j=1}^{N_0}J_s(w_j)^2
\le C\sum_{j=1}^{N_0}\int_{\mathbb{R}^n_+}
t^{-2\gamma}|\nabla w_j(x',t)|^2\,dx'\,dt
\le C\int_\Omega\delta^{-2\gamma}
\big(|\nabla v|^2+|v|^2\big)\,dx
\le C\|v\|_{H^{3s/2}(\Omega)}^2.
\end{equation*}
This finishes the proof of Lemma~\ref{l10}.
\end{proof}

Let $u\in H^{3s/2}(\Omega)$. For each $j\in\{1,\ldots,N_0\}$ and any
$z\in U_j^1\cap\Omega^c$, define
\begin{equation*}
m_j(z):=u(z)-u(R_jz).
\end{equation*}

\begin{lemma}\label{l11}
Let $n\ge2$, let $s\in(0,1)$, and let $\Omega\subset\mathbb{R}^n$ be a bounded
Lipschitz domain. Assume that $q\in L^2(\Omega^c,\sigma_s)$ satisfies
$\int_{\Omega^c}q\,d\sigma_s=0$, $g:=a_sq$, and
$u\in\mathcal V$ is the weak solution in the sense of
Definition~\ref{def:weak}. Then
\begin{equation*}
\sum_{j=1}^{N_0}\int_{U_j^1\cap\Omega^c}
\delta^{-3s}|m_j|^2\,dz
\le C\left[\|q\|_{L^2(\Omega^c,\sigma_s)}^2
+\|u\|_{H^{3s/2}(\Omega)}^2\right],
\end{equation*}
where the positive constant $C$ depends only on $n$, $s$, and $\Omega$.
\end{lemma}

\begin{proof}
By \eqref{e2.25}, we find that
\begin{equation}\label{e4.31}
m_j=\frac{a_sq}{c_{(n,s)}\rho}+Bu-u\circ R_j
\end{equation}
on $U_j^1\cap\Omega^c$. From the facts that $a_s\sim\delta^{-s}$ and
$\rho\sim\delta^{-2s}$, it follows that, for any
$z\in U_j^1\cap\Omega^c$,
\begin{equation*}
\delta(z)^{-3s}\left|\frac{a_s(z)q(z)}{c_{(n,s)}\rho(z)}\right|^2
\le C\delta(z)^{-s}|q(z)|^2
\le Ca_s(z)|q(z)|^2.
\end{equation*}
By integrating over the charts, we conclude that
\begin{equation}\label{e4.32}
\int_{U_j^1\cap\Omega^c}\delta(z)^{-3s}
\left|\frac{a_s(z)q(z)}{c_{(n,s)}\rho(z)}\right|^2\,dz
\le C\|q\|_{L^2(\Omega^c,\sigma_s)}^2.
\end{equation}
From the definition $\rho(z)=\int_\Omega K(z,x)\,dx$, we deduce that, for
any $z\in U_j^1\cap\Omega^c$,
\begin{equation*}
Bu(z)-u(R_jz)=\rho(z)^{-1}\int_\Omega K(z,x)[u(x)-u(R_jz)]\,dx.
\end{equation*}
By $\rho(z)^{-1}\le C\delta(z)^{2s}$ in Lemma~\ref{l2}(i), we find that
\begin{equation*}
\delta(z)^{-3s}\left|Bu(z)-u(R_jz)\right|^2
\le C\delta(z)^{-s}\int_\Omega K(z,x)\left|u(x)-u(R_jz)\right|^2\,dx.
\end{equation*}

Now, let $x\in U_j^2\cap\Omega$ and $z\in U_j^1\cap\Omega^c$. Write
$x=F_j(X)$, $z=F_j(y',-r)$, and $R_jz=F_j(Y)$, where $X:=(x',t)$ and
$Y:=(y',r)$ with $t,r\in(0,\infty)$. From \eqref{e4.25}, \eqref{e4.26}, and
\eqref{e4.27}, it follows that
\begin{equation*}
u(x)=w_j(X),\quad u(R_jz)=w_j(Y),\quad
\delta(z)^{-s}\le Cr^{-s},\quad\text{and}\quad
K(z,x)\le CK_{\rm ref}(X,Y).
\end{equation*}
By this and Lemma~\ref{l10}, we conclude that
\begin{align}\label{e4.33}
&\int_{U_j^1\cap\Omega^c}\delta(z)^{-s}
\int_{U_j^2\cap\Omega}K(z,x)\left|u(x)-u(R_jz)\right|^2\,dx\,dz\nonumber\\
&\quad\leq C\iint_{\mathbb{R}^n_+\times\mathbb{R}^n_+}r^{-s}K_{\rm ref}(X,Y)
\left|w_j(X)-w_j(Y)\right|^2\,dX\,dY=CJ_s(w_j)^2\leq C\|u\|_{H^{3s/2}(\Omega)}^2.
\end{align}
From \eqref{e2.4} with $t:=s/2$ and the fact that
$\operatorname{dist}(U_j^1,\mathbb{R}^n\setminus U_j^2)>0$, we deduce
that
\begin{align}\label{e4.34}
&\int_{U_j^1\cap\Omega^c}\delta(z)^{-s}
\int_{\Omega\setminus U_j^2}K(z,x)\left|u(x)-u(R_jz)\right|^2\,dx\,dz\nonumber\\
&\quad\leq C\|u\|_{L^2(\Omega)}^2
+C\int_{U_j^1\cap\Omega^c}\delta(z)^{-s}\left|u(R_jz)\right|^2\,dz\nonumber\\
&\quad\leq C\|u\|_{L^2(\Omega)}^2
+C\int_\Omega\delta^{-s}|u|^2\,dx\leq C\|u\|_{H^{s/2}(\Omega)}^2.
\end{align}
Then, by \eqref{e4.31}, \eqref{e4.32}, \eqref{e4.33}, \eqref{e4.34}, and
Lemma~\ref{l10}, we find that
\begin{equation*}
\begin{aligned}
\sum_{j=1}^{N_0}\int_{U_j^1\cap\Omega^c}\delta^{-3s}|m_j|^2\,dz
&\le C\left[\|q\|_{L^2(\Omega^c,\sigma_s)}^2
+\sum_{j=1}^{N_0}J_s(w_j)^2+\|u\|_{H^{s/2}(\Omega)}^2\right]\\
&\le C\left[\|q\|_{L^2(\Omega^c,\sigma_s)}^2
+\|u\|_{H^{3s/2}(\Omega)}^2\right].
\end{aligned}
\end{equation*}
This finishes the proof of Lemma~\ref{l11}.
\end{proof}

\section{Weighted gradient and non-tangential maximal estimates}\label{s5}

In this section, we first establish an estimate for tails and reconstruct the
derivatives for the weak solution to \eqref{e1.1}, and then finish the proofs of Theorems \ref{t1} and \ref{t2}.
Throughout this section, assume that $n\geq2$, $s\in(0,1)$, and
$\Omega\subset\mathbb{R}^n$ is a bounded Lipschitz domain. For
$q\in L^2(\Omega^c,\sigma_s)$ and $u\in\mathcal{V}$ with
$u\in H^{3s/2}(\Omega)$, define
\begin{equation*}
\mathcal A_s(u,q):=\left[\|u\|_{H^{3s/2}(\Omega)}^2
+\|q\|_{L^2(\Omega^c,\sigma_s)}^2\right]^{1/2}.
\end{equation*}
From \eqref{e2.19}, it follows that
\begin{equation}\label{e5.2}
\mathcal{E}_\Omega(u,u)+\|u\|_{L^2(\Omega^c,\sigma_s)}^2
\leq C\|q\|_{L^2(\Omega^c,\sigma_s)}^2
\leq C\mathcal A_s(u,q)^2.
\end{equation}
By the interior Sobolev estimate \eqref{e3.34}, we conclude that
$\mathcal A_s(u,q)\leq C\|q\|_{L^2(\Omega^c,\sigma_s)}$.
Whenever $w_j$ is used in this section, it is constructed by
\eqref{e4.27} from the values of $u$ in $\Omega$.

\subsection{Tail difference}\label{s5.1}

Fix $\kappa\in(0,1/4)$ and define $\mathcal T_\kappa(u)$ by setting
\begin{equation}\label{e5.3}
\mathcal T_\kappa(u):=\int_\Omega\delta(x)^{-s}
\int_{\{y\in\mathbb{R}^n:|y-x|>\kappa\delta(x)\}}K(x,y)
|u(y)-u(x)|^2\,dy\,dx.
\end{equation}

\begin{lemma}\label{l12}
Let $n\ge2$, let $s\in(0,1)$, and let $\Omega\subset\mathbb{R}^n$ be a bounded
Lipschitz domain. Assume that $q\in L^2(\Omega^c,\sigma_s)$ satisfies
$\int_{\Omega^c}q\,d\sigma_s=0$, let $g:=a_sq$, and let
$u\in\mathcal V$ be the weak solution in the sense of
Definition~\ref{def:weak}. Then, for any $\kappa\in(0,1/4)$,
\begin{equation*}
\mathcal T_\kappa(u)\le C_\kappa\mathcal A_s(u,q)^2,
\end{equation*}
where the positive constant $C_\kappa$ depends only on $n$, $s$, $\kappa$,
and $\Omega$.
\end{lemma}

\begin{proof}
For any $x\in U_j^0\cap\Omega$, we split the integration in $y$ into the
three regions $U_j^1\cap\Omega$, $U_j^1\cap\Omega^c$, and
$\mathbb{R}^n\setminus U_j^1$.

We first consider the case $y\in U_j^1\cap\Omega$. Let $x:=F_j(x',t)$ and
$y:=F_j(y',r)$ with $t,r\in(0,\infty)$. From the facts that $|x-y|>\kappa\delta(x)$ and
$t\sim\delta(x)$, we deduce that $t\le C_\kappa|x-y|$. By \eqref{e4.25}, we
find that
\begin{equation}\label{e5.4}
|t-r|+|x'-y'|\le C|x-y|.
\end{equation}
From this and the fact that $t+r\le2t+|t-r|$, it follows that
\begin{equation}\label{e5.5}
t+r\le C_\kappa|x-y|.
\end{equation}
By \eqref{e5.4} and \eqref{e5.5}, we conclude that
$K(x,y)\le C_\kappa K_{\rm ref}(X,Y)$, where $X:=(x',t)$ and $Y:=(y',r)$.
From the assumption that $\psi_j=1$ on $U_j^2$, we deduce that
\begin{align}\label{e5.6}
&\int_{U_j^0\cap\Omega}\delta(x)^{-s}\int_{\substack{y\in U_j^1\cap\Omega\\|y-x|>\kappa\delta(x)}}
K(x,y)|u(y)-u(x)|^2\,dy\,dx\nonumber\\
&\quad\le C_\kappa\iint_{\mathbb{R}^n_+\times\mathbb{R}^n_+}t^{-s}K_{\rm ref}(X,Y)
|w_j(X)-w_j(Y)|^2\,dX\,dY=C_\kappa J_s(w_j)^2.
\end{align}

We next consider the case $y\in U_j^1\cap\Omega^c$. By the definition of
$m_j$ in Lemma~\ref{l11}, we find that, for any
$y=F_j(y',-r)\in U_j^1\cap\Omega^c$ with $r\in(0,\infty)$,
\begin{equation}\label{e5.7}
u(y)-u(x)=m_j(y)+u(R_jy)-u(x).
\end{equation}
From \eqref{e4.25}, it follows that
\begin{equation}\label{equ5.8}
\int_{U_j^0\cap\Omega}\delta(x)^{-s}
\int_{U_j^1\cap\Omega^c}K(x,y)|u(R_jy)-u(x)|^2\,dy\,dx
\le CJ_s(w_j)^2.
\end{equation}
By Fubini's theorem, we conclude that
\begin{equation}\label{e5.8}
\int_{U_j^0\cap\Omega}\delta(x)^{-s}
\int_{U_j^1\cap\Omega^c}K(x,y)|m_j(y)|^2\,dy\,dx
\le\int_{U_j^1\cap\Omega^c}\delta(y)^{-3s}|m_j(y)|^2\,dy.
\end{equation}
Then, from \eqref{e5.7}, \eqref{equ5.8}, \eqref{e5.8}, and the
inequality $|A+B|^2\le2|A|^2+2|B|^2$, we deduce that
\begin{align}\label{e5.9}
&\int_{U_j^0\cap\Omega}\delta(x)^{-s}
\int_{U_j^1\cap\Omega^c}K(x,y)|u(y)-u(x)|^2\,dy\,dx\notag\\
&\quad\le C\left(J_s(w_j)^2
+\int_{U_j^1\cap\Omega^c}\delta^{-3s}|m_j|^2\,dy\right).
\end{align}

We finally consider the case $y\in\mathbb{R}^n\setminus U_j^1$. By the
positive lower bound of $a_s$ on bounded exterior sets and the asymptotic
relation $K(x,y)\sim a_s(y)$ at infinity, we find that, for any
$x\in U_j^0\cap\Omega$ and $y\in\Omega^c\setminus U_j^1$,
\begin{equation}\label{e5.10}
K(x,y)\le C_j a_s(y).
\end{equation}
Moreover, from Hardy's inequality, it follows that
\begin{equation}\label{e5.12}
\int_\Omega\delta^{-s}|u|^2\,dx
\le C\|u\|_{H^{s/2}(\Omega)}^2.
\end{equation}
Then, by \eqref{e5.10}, \eqref{e5.12}, \eqref{e5.2}, and the fact that
$\int_\Omega\delta(x)^{-s}\,dx<\infty$, we conclude that
\begin{align}\label{e5.13}
&\int_{U_j^0\cap\Omega}\delta(x)^{-s}
\int_{\mathbb{R}^n\setminus U_j^1}K(x,y)|u(y)-u(x)|^2\,dy\,dx\notag\\
&\quad\le C_j\left(\|u\|_{L^2(\Omega)}^2
+\|u\|_{L^2(\Omega^c,\sigma_s)}^2
+\int_\Omega\delta(x)^{-s}|u(x)|^2\,dx\right)
\le C_j\mathcal A_s(u,q)^2.
\end{align}

From the fact that there exists a constant $d_0>0$ such that
$\delta(x)\ge d_0$ for any
$x\in\Omega\setminus(\bigcup_{j=1}^{N_0}U_j^0)$, we deduce that
\begin{align}\label{e5.14}
&\int_{\Omega\setminus\bigcup_{j=1}^{N_0}U_j^0}\delta(x)^{-s}
\int_{\{y\in\mathbb{R}^n:|y-x|>\kappa\delta(x)\}}
K(x,y)|u(y)-u(x)|^2\,dy\,dx\notag\\
&\quad\le d_0^{-s}\int_\Omega\int_{\mathbb{R}^n}
K(x,y)|u(y)-u(x)|^2\,dy\,dx
\le\frac{2d_0^{-s}}{c_{(n,s)}}\mathcal E_\Omega(u,u).
\end{align}
Then, by \eqref{e5.6}, \eqref{e5.9}, \eqref{e5.13}, and \eqref{e5.14},
together with Lemmas~\ref{l10} and~\ref{l11}, we find that
\begin{align*}
\mathcal T_\kappa(u)
&\le C_\kappa\Bigg[\sum_{j=1}^{N_0}J_s(w_j)^2
+\sum_{j=1}^{N_0}\int_{U_j^1\cap\Omega^c}
\delta^{-3s}|m_j|^2\,dz\\
&\quad+\|u\|_{H^{s/2}(\Omega)}^2+\mathcal E_\Omega(u,u)
+\|u\|_{L^2(\Omega^c,\sigma_s)}^2\Bigg]
\le C_\kappa\mathcal A_s(u,q)^2.
\end{align*}
This finishes the proof of Lemma~\ref{l12}.
\end{proof}

\subsection{Weighted gradient estimate}\label{s5.2}

Let $w_s(y):=(1+|y|)^{-n-2s}$ for any $y\in\mathbb{R}^n$ and
\begin{equation*}
L^1(\mathbb{R}^n,w_s\,dy):=\left\{f:\mathbb{R}^n\to\mathbb{R}\text{ measurable}:\
\int_{\mathbb{R}^n}|f(y)|w_s(y)\,dy<\infty\right\}.
\end{equation*}
Assume that $u\in\mathcal V$ is the weak solution to problem \eqref{e1.1} in the sense of
Definition~\ref{def:weak}. Then, by H\"older's inequality and \eqref{e2.19},
we conclude that
\begin{equation*}
\int_{\Omega^c}|u|w_s\,dy
\le\|u\|_{L^2(\Omega^c,\sigma_s)}
\left(\int_{\Omega^c}\frac{w_s^2}{a_s}\,dy\right)^{1/2}
<\infty.
\end{equation*}
From the facts that $a_s^{-1}\le C\delta^s$ near the boundary, the decay
$w_s^2/a_s\sim|y|^{-n-2s}$ at infinity, and $u\in L^2(\Omega)$, we deduce
that $u\in L^1(\mathbb{R}^n,w_s\,dy)$. By the definition of the weak
solution, symmetry, and Fubini's theorem, we find that, for any
$\phi\in C_{\rm c}^\infty(\Omega)$,
\begin{equation}\label{equ5.15}
\int_{\mathbb{R}^n}u(-\Delta)^s\phi\,dx=0.
\end{equation}

The following lemma is a variant of the Cauchy estimate in
\cite[Lemma~3.1]{Fa16}.

\begin{lemma}\label{l13}
Let $n\ge2$, $s\in(0,1)$, and let $\Omega\subset\mathbb{R}^n$ be open.
Assume that $u\in L^1(\mathbb{R}^n,w_s\,dy)$ satisfies \eqref{equ5.15}.
Then, for any $m\in\mathbb N$ and any ball $B(x_0,2\varrho)\Subset\Omega$,
\begin{equation*}
\varrho^{2m}\sup_{x\in B(x_0,\varrho/2)}|D^mu(x)|^2
\le C_{n,s,m}\varrho^{2s}
\int_{|y-x_0|>\varrho}K(x_0,y)|u(y)-u(x_0)|^2\,dy,
\end{equation*}
where the positive constant $C_{n,s,m}$ depends only on $n$, $s$, and $m$,
and is independent of $\Omega$, $u$, $x_0$, and $\varrho$.
\end{lemma}

\begin{proof}
Choose a nonnegative function $\eta\in C_{\rm c}^\infty(B({\bf0},1))$ such that
$\int_{\mathbb{R}^n}\eta\,dx=1$, and define
$\eta_\varepsilon(y):=\varepsilon^{-n}\eta(y/\varepsilon)$ and
$u_\varepsilon:=u*\eta_\varepsilon$. Then
$u_\varepsilon\in C(\mathbb{R}^n)\cap L^1(\mathbb{R}^n,w_s\,dy)$ and
\begin{equation}\label{equ5.17}
\lim_{\varepsilon\to0^+}\|u_\varepsilon-u\|_{L^1(\mathbb{R}^n,w_s\,dy)}=0.
\end{equation}
Fix $\varepsilon\in(0,\operatorname{dist\,}(B(x_0,2\varrho),\Omega^c))$ and
let $\check\eta_\varepsilon(y):=\eta_\varepsilon(-y)$. From Fubini's
theorem and \eqref{equ5.15}, it follows that, for any
$\phi\in C_{\rm c}^\infty(B(x_0,2\varrho))$,
\begin{equation}\label{e5.17}
\langle(-\Delta)^su_\varepsilon,\phi\rangle
=\int_{\mathbb{R}^n}u(y)(-\Delta)^s
(\check\eta_\varepsilon*\phi)(y)\,dy=0.
\end{equation}
Here $\check\eta_\varepsilon*\phi\in C_{\rm c}^\infty(\Omega)$ and
$|(-\Delta)^s(\check\eta_\varepsilon*\phi)|\le Cw_s$. By \eqref{e5.17},
we conclude that $(-\Delta)^su_\varepsilon=0$ in $B(x_0,2\varrho)$. From
the ball Poisson formula (see, for instance,
\cite[Definition~1.8 and Theorem~2.10]{Bu16}), we deduce that, for any
$x\in B(x_0,\varrho/2)$ and any $R\in(\varrho,2\varrho)$,
\begin{equation}\label{equ5.19}
u_\varepsilon(x)=\int_{\mathbb{R}^n}P_R(x,y)u_\varepsilon(y)\,dy,
\end{equation}
where $P_R(x,y):=0$ for any $y\in B(x_0,R)$, and
\begin{equation*}
P_R(x,y):=p_{n,s}\frac{(R^2-|x-x_0|^2)^s}
{(|y-x_0|^2-R^2)^s|x-y|^n}
\end{equation*}
for any $y\in\mathbb{R}^n\setminus B(x_0,R)$, with
$p_{n,s}:=\Gamma(n/2)\sin(\pi s)/\pi^{n/2+1}$. By \eqref{equ5.19}, we find
that
\begin{equation}\label{equ5.20}
u_\varepsilon(x)=\frac{1}{\varrho}\int_{\mathbb{R}^n}
\int_{\varrho}^{\min\{|y-x_0|,2\varrho\}}P_R(x,y)u_\varepsilon(y)
\,dR\,dy.
\end{equation}

For any $x\in B(x_0,\varrho/2)$,
$y\in\mathbb{R}^n\setminus B(x_0,\varrho)$, and
$R\in(\varrho,2\varrho)$, we have
\begin{equation*}
\frac{3}{4}\varrho^2\le R^2-|x-x_0|^2\le4\varrho^2
\quad\text{and}\quad
|x-y|\ge\frac{1}{2}|y-x_0|,
\end{equation*}
which further imply that, for any
$j\in\{0,\ldots,m\}$,
\begin{equation*}
\left|D_x^j(R^2-|x-x_0|^2)^s\right|\le C_j\varrho^{2s-j}
\quad\text{and}\quad
\left|D_x^{m-j}\left(|x-y|^{-n}\right)\right|\le C_m\varrho^{j-m}
|y-x_0|^{-n}.
\end{equation*}
By this and the definition of $P_R(x,y)$, we conclude that, for any multi-index $\alpha$,
$x\in B(x_0,\varrho/2)$, and $y\in\mathbb{R}^n$,
\begin{equation}\label{e5.23}
\left|\int_{\varrho}^{\min\{|y-x_0|,2\varrho\}}\partial_x^\alpha
P_R(x,y)\,dR\right|\le C\varrho^{2s-m}|y-x_0|^{-n-2s}
\mathbf{1}_{\{|y-x_0|>\varrho\}}.
\end{equation}
Then, from \eqref{e5.23} and \eqref{equ5.20}, we deduce that, for any
multi-index $\alpha$,
\begin{align*}
\sup_{x\in B(x_0,\varrho/2)}\left|\partial^\alpha u_\varepsilon(x)
-\int_{\mathbb{R}^n}\int_{\varrho}^{\min\{|y-x_0|,2\varrho\}}
\partial_x^\alpha P_R(x,y)u(y)\,dR\,dy\right|\le C_{\varrho,\alpha}\|u_\varepsilon-u\|_{L^1(\mathbb{R}^n,w_s\,dy)},
\end{align*}
which, together with \eqref{equ5.17} and \eqref{equ5.20}, further implies that, for
any multi-index $\alpha$,
\begin{equation}\label{equ5.26}
\partial^\alpha u(x)=\frac{1}{\varrho}\int_{\mathbb{R}^n}
\int_{\varrho}^{\min\{|y-x_0|,2\varrho\}}\partial^\alpha_xP_R(x,y)u(y)
\,dR\,dy.
\end{equation}

Then, by \eqref{equ5.26}, we conclude that, for any $m\in\mathbb N$ and any
$x\in B(x_0,\varrho/2)$,
\begin{equation*}
\int_{\mathbb{R}^n}\int_{\varrho}^{\min\{|y-x_0|,2\varrho\}}
D^m_xP_R(x,y)\,dR\,dy=0
\end{equation*}
and
\begin{equation*}
D^mu(x)=\frac{1}{\varrho}\int_{\mathbb{R}^n}
\int_{\varrho}^{\min\{|y-x_0|,2\varrho\}}D^m_xP_R(x,y)u(y)\,dR\,dy.
\end{equation*}
From this and \eqref{e5.23}, we infer that
\begin{equation*}
|D^mu(x)|\le C_m\varrho^{2s-m}\int_{|y-x_0|>\varrho}
K(x_0,y)|u(y)-u(x_0)|\,dy,
\end{equation*}
which, combined with H\"older's inequality and the fact that
$\int_{|y-x_0|>\varrho}K(x_0,y)\,dy=C\varrho^{-2s}$, implies that
\begin{equation*}
\varrho^{2m}\sup_{x\in B(x_0,\varrho/2)}|D^mu(x)|^2
\le C_m\varrho^{2s}\int_{|y-x_0|>\varrho}
K(x_0,y)|u(y)-u(x_0)|^2\,dy.
\end{equation*}
This finishes the proof of Lemma~\ref{l13}.
\end{proof}

\begin{proposition}\label{t5}
Let $n\ge2$, let $s\in(0,1)$, and let $\Omega\subset\mathbb{R}^n$ be a bounded
Lipschitz domain. Assume that $q\in L^2(\Omega^c,\sigma_s)$ satisfies
$\int_{\Omega^c}q\,d\sigma_s=0$, $g:=a_sq$, and $u\in\mathcal V$ is
the weak solution to problem \eqref{e1.1} in the sense of Definition~\ref{def:weak}. Then, for any
$m\in\mathbb N$,
\begin{equation*}
\int_\Omega\delta^{2m-3s}|D^mu|^2\,dx
\le C_m\mathcal T_{1/8}(u),
\end{equation*}
where the positive constant $C_m$ depends only on $n$, $s$, $m$, and
$\Omega$.
\end{proposition}

\begin{proof}
Applying Lemma~\ref{l13} with $\varrho:=\delta(x)/8$, we conclude that
\begin{equation*}
\int_\Omega\delta^{2m-3s}|D^mu|^2\,dx
\le C_m\int_\Omega\delta(x)^{-s}
\int_{|y-x|>\delta(x)/8}K(x,y)|u(y)-u(x)|^2\,dy\,dx
=C_m\mathcal T_{1/8}(u).
\end{equation*}
This finishes the proof of Proposition~\ref{t5}.
\end{proof}

\subsection{Non-tangential maximal estimate}\label{s5.3}

\begin{proposition}\label{l14}
Let $n\ge2$, let $s\in(0,1)$, and let $\Omega\subset\mathbb{R}^n$ be a bounded
Lipschitz domain. Assume that $q\in L^2(\Omega^c,\sigma_s)$ satisfies
$\int_{\Omega^c}q\,d\sigma_s=0$, $g:=a_sq$, and 
$u\in\mathcal V$ is the weak solution to problem \eqref{e1.1} in the sense of
Definition~\ref{def:weak}. Then
\begin{equation*}
\left\|\left(\delta^{1-s}|\nabla u|\right)^*\right\|_{L^2(\Omega^c,\sigma_s)}^2
\le C\,\mathcal T_{1/8}(u),
\end{equation*}
where the positive constant $C$ depends only on $n$, $s$, and $\Omega$.
\end{proposition}

\begin{proof}
By the Whitney decomposition theorem (see, for instance,
\cite[Chapter~VI, Section~1]{Stein70}), we find that there exists a family
$\mathcal W$ of maximal dyadic cubes $Q\subset\Omega$ satisfying
$\delta(Q):=\inf_{x\in Q}\delta(x)\ge8d_Q,$
where $d_Q:=\operatorname{diam\,}(Q)=\sqrt{n}\,\ell_Q$ and $\ell_Q$ denotes
the side length of $Q$. That is, $\Omega=\bigcup_{Q\in\mathcal W}Q$ and,
for any distinct $Q,Q'\in\mathcal W$,
$\operatorname{int}Q\cap\operatorname{int}Q'=\varnothing$. Then, for any
$x\in Q$,
\begin{equation}\label{e5.25}
8d_Q\le\delta(x)<10d_Q.
\end{equation}
For any cube $Q\in\mathcal W$, let
\begin{equation*}
S_Q:=\left\{z\in\Omega^c:\Gamma(z)\cap Q\ne\varnothing\right\}.
\end{equation*}
From \eqref{e5.25}, we deduce that, for the center $x_Q$ of $Q$ and any
$z\in S_Q$, one has $|z-x_Q|\le C\ell_Q$ and $\delta(z)\le C\ell_Q$,
which, combined with the definition of $\sigma_s$, implies that, for any cube $Q$ with $C\ell_Q<r_0$,
\begin{equation}\label{e5.26}
\sigma_s(S_Q)
\le C\ell_Q^{n-1}\int_0^{C\ell_Q}r^{-s}\,dr
\le C\ell_Q^{n-s}.
\end{equation}
Furthermore, by the fact that $\sigma_s(\Omega^c)<\infty$, we conclude that
\eqref{e5.26} also holds for any cube $Q$ with $C\ell_Q\ge r_0$.

From the definition of $S_Q$ and the fact that $\delta(x)\sim\ell_Q$ for any
$x\in Q$, it follows that, for any $z\in\Omega^c$,
\begin{equation*}
\left[\left(\delta^{1-s}|\nabla u|\right)^*(z)\right]^2
\le C\sum_{Q\in\mathcal W}\mathbf{1}_{S_Q}(z)\ell_Q^{2-2s}
\sup_{x\in Q}|\nabla u(x)|^2,
\end{equation*}
which further implies that
\begin{equation}\label{e5.26x}
\left\|\left(\delta^{1-s}|\nabla u|\right)^*\right\|_{L^2(\Omega^c,\sigma_s)}^2
\le C\sum_{Q\in\mathcal W}\ell_Q^{n+2-3s}\sup_{x\in Q}|\nabla u(x)|^2.
\end{equation}

For any cube $Q$ and any $x\in Q$, let
$M_Q:=\sup_{z\in Q}|\nabla u(z)|^2$
and
\begin{equation*}
I_Q(x):=\int_{\{y\in\mathbb{R}^n:|y-x|>d_Q\}}
K(x,y)|u(y)-u(x)|^2\,dy.
\end{equation*}
By Lemma~\ref{l13} with $m:=1$ and the fact that $|z-x|\le d_Q$ for any
$z,x\in Q$, we find that, for any $x\in Q$,
\begin{equation*}
M_Q\le Cd_Q^{2s-2}I_Q(x),
\end{equation*}
which, combined with \eqref{e5.25}, further implies that
\begin{equation*}
\ell_Q^{n+2-3s}M_Q
\le C\int_Q\delta(x)^{-s}
\int_{\{y\in\mathbb{R}^n:|y-x|>\delta(x)/10\}}
K(x,y)|u(y)-u(x)|^2\,dy\,dx.
\end{equation*}
From this and \eqref{e5.26x}, we deduce that
\begin{align*}
\left\|\left(\delta^{1-s}|\nabla u|\right)^*\right\|_{L^2(\Omega^c,\sigma_s)}^2
&\le C\sum_{Q\in\mathcal W}\ell_Q^{n+2-3s}M_Q\\
&\le C\sum_{Q\in\mathcal W}\int_Q\delta(x)^{-s}
\int_{|y-x|>\delta(x)/8}K(x,y)|u(y)-u(x)|^2\,dy\,dx\\
&=C\mathcal T_{1/8}(u).
\end{align*}
This finishes the proof of Proposition~\ref{l14}.
\end{proof}

\subsection{Proofs of the main results}\label{s5.4}

\begin{proof}[Proof of Theorem~\ref{t1}]
From Lemma~\ref{l3}\textup{(ii)}, it follows that there exists a unique
solution $u$ to problem \eqref{e1.1}. Let $q:=g/a_s$ and define, for any
$w\in H^{s/2}(\Omega)$,
\begin{equation*}
f(w):=F_q(w)+(u,w)_{L^2(\Omega)}.
\end{equation*}
By Lemmas~\ref{l3}\textup{(ii)} and~\ref{l4}\textup{(iv)},
\eqref{e3.34}, and the definition of $\mathcal A_s$, we conclude that
\begin{equation*}
\|f\|_{Y_{s/2}}\le C\left[\|g\|_{\mathcal X_s(\Omega^c)}+\|u\|_{L^2(\Omega)}\right]
\le C\|g\|_{\mathcal X_s(\Omega^c)},
\end{equation*}
\begin{equation*}
\|u\|_{H^{3s/2}(\Omega)}\le C\|g\|_{\mathcal X_s(\Omega^c)},
\quad\text{and}\quad
\mathcal A_s(u,q)\le C\|g\|_{\mathcal X_s(\Omega^c)}.
\end{equation*}
This proves \eqref{e1.12x}.

From Lemma~\ref{l12} and Propositions~\ref{t5} and ~\ref{l14}, we deduce that
\begin{equation*}
\int_\Omega\delta^{2-3s}|\nabla u|^2\,dx
+\left\|\left(\delta^{1-s}|\nabla u|\right)^*\right\|_{L^2(\Omega^c,\sigma_s)}^2
\le C\mathcal T_{1/8}(u)\le C\|g\|_{\mathcal X_s(\Omega^c)}^2.
\end{equation*}
This proves \eqref{e1.13} and \eqref{e1.14}, and hence completes the proof
of Theorem~\ref{t1}.
\end{proof}

\begin{proof}[Proof of Theorem~\ref{t2}]
By Theorem~\ref{t1}, Proposition~\ref{t5}, and \eqref{e1.4},
we conclude that, for any $m\in\mathbb N$,
\begin{equation*}
\int_\Omega\delta^{2m-3s}|D^mu|^2\,dx
\le C_m\mathcal T_{1/8}(u)
\le C_m\mathcal A_s(u,q)^2\le C_m\|q\|_{L^2(\Omega^c,\sigma_s)}^2
=C_m\|g\|_{\mathcal X_s(\Omega^c)}^2.
\end{equation*}
This finishes the proof of Theorem~\ref{t2}.
\end{proof}

\medskip

\noindent\textbf{Data Availability}\quad
This paper has no associated data.

\section*{Declarations}

\noindent\textbf{Conflict of interest}\quad The authors declare
that they have no conflict of interest.

\medskip

\noindent\textbf{Declaration of AI Use}\quad
The authors developed the initial theoretical framework, and ChatGPT was used to assist in 
discussing proof strategies, verifying mathematical arguments, and preparing the
manuscript. The final manuscript was written and verified by the authors, 
who take full responsibility for its content.

\bigskip

\noindent Wenxian Ma and Sibei Yang

\medskip

\noindent School of Mathematics and Statistics, Gansu Key Laboratory of
Applied Mathematics and Complex Systems, Lanzhou University,
Lanzhou 730000, People's Republic of China

\smallskip

\noindent {\it E-mails} \texttt{mawenxian@lzu.edu.cn} (W. Ma)

\hspace{0.888cm}\texttt{yangsb@lzu.edu.cn} (S. Yang)

\begin{thebibliography}{99}

\bibitem{AFR23}
A. Audrito, J.-C. Felipe-Navarro and X. Ros-Oton,
The Neumann problem for the fractional Laplacian: regularity up to
the boundary,
Ann. Sc. Norm. Super. Pisa Cl. Sci. (5) 24 (2023), 1155--1222.

\vspace{-0.3cm}

\bibitem{BL76}
J. Bergh and J. L\"ofstr\"om,
Interpolation Spaces: An Introduction,
Grundlehren der mathematischen Wissenschaften, Vol. 223,
Springer-Verlag, Berlin--Heidelberg, 1976.

\vspace{-0.3cm}

\bibitem{BB22}
F. M. Bersetche and J. P. Borthagaray,
Finite element approximation of fractional Neumann problems,
IMA J. Numer. Anal. 42 (2022), 3207--3240.

\vspace{-0.3cm}

\bibitem{BLN23}
J. P. Borthagaray, W. Li and R. H. Nochetto,
Quasi-linear fractional-order operators in Lipschitz domains,
SIAM J. Math. Anal. 56 (2024), 4006--4039.

\vspace{-0.3cm}

\bibitem{BN23}
J. P. Borthagaray and R. H. Nochetto,
Besov regularity for the Dirichlet integral fractional Laplacian
in Lipschitz domains,
J. Funct. Anal. 284 (2023), Article 109829, 33 pp.

\vspace{-0.3cm}

\bibitem{Bu16}
C. Bucur,
Some observations on the Green function for the ball in the
fractional Laplace framework,
Commun. Pure Appl. Anal. 15 (2016), 657--699.

\vspace{-0.3cm}

\bibitem{cs07} L. Caffarelli and L. Silvestre, An extension problem related to the fractional Laplacian,
Comm. Partial Differential Equations 32 (2007), 1245--1260.

\vspace{-0.3cm}

\bibitem{CFR26}
R. Colombo, X. Fern\'andez-Real and X. Ros-Oton,
The fractional Laplacian in Lipschitz domains: Dahlberg's theorem
and $L^2$-solvability,
arXiv:2607.14909.

\vspace{-0.3cm}

\bibitem{Dah77}
B. E. Dahlberg,
Estimates of harmonic measure,
Arch. Rational Mech. Anal. 65 (1977), 275--288.

\vspace{-0.3cm}

\bibitem{DK87}
B. E. Dahlberg and C. E. Kenig,
Hardy spaces and the Neumann problem in $L^p$ for Laplace's equation
in Lipschitz domains,
Ann. of Math. (2) 125 (1987), 437--465.

\vspace{-0.3cm}

\bibitem{DRV17}
S. Dipierro, X. Ros-Oton and E. Valdinoci,
Nonlocal problems with Neumann boundary conditions,
Rev. Mat. Iberoam. 33 (2017), 377--416.

\vspace{-0.3cm}

\bibitem{DRVW25}
S. Dipierro, X. Ros-Oton, E. Valdinoci and M. Weidner,
The Neumann problem for the fractional Laplacian: optimal regularity
via the Mellin transform, arXiv:2510.13340v2,
to appear in Comm. Amer. Math. Soc.

\vspace{-0.3cm}

\bibitem{DNPV12}
E. Di Nezza, G. Palatucci and E. Valdinoci,
Hitchhiker's guide to the fractional Sobolev spaces,
Bull. Sci. Math. 136 (2012), 521--573.

\vspace{-0.3cm}

\bibitem{dglz12} Q. Du, M. Gunzburger, R. B. Lehoucq and K. Zhou, Analysis and approximation of nonlocal diffusion
problems with volume constraints, SIAM Review 54 (2012), 667--696.

\vspace{-0.3cm}

\bibitem{fjr78}E. B. Fabes, M. Jr. Jodeit and N. M. Rivi\'ere, Potential techniques for
boundary value problems on  $C^1$-domains, Acta Math. 141 (1978), 165--186.

\vspace{-0.3cm}

\bibitem{Fa16}
M. M. Fall,
Entire $s$-harmonic functions are affine,
Proc. Amer. Math. Soc. 144 (2016), 2587--2592.

\vspace{-0.3cm}

\bibitem{FK24}
G. Foghem and M. Kassmann,
A general framework for nonlocal Neumann problems,
Commun. Math. Sci. 22 (2024), 15--66.

\vspace{-0.3cm}

\bibitem{fls16} R. L. Frank, E. Lenzmann and L. Silvestre, Uniqueness of radial solutions for the fractional Laplacian,
Comm. Pure Appl. Math. 69 (2016), 1671--1726.

\vspace{-0.3cm}

\bibitem{FS10}
R. L. Frank and R. Seiringer,
Sharp fractional Hardy inequalities in half-spaces,
in Around the Research of Vladimir Maz'ya I,
International Mathematical Series, Vol. 11,
Springer, New York, 2010, 161--167.

\vspace{-0.3cm}

\bibitem{GH24}
F. Grube and T. Hensiek,
Robust nonlocal trace spaces and Neumann problems,
Nonlinear Anal. 241 (2024), Paper No. 113481, 35 pp.

\vspace{-0.3cm}

\bibitem{jsw20} S. Jarohs, A. Salda$\mathrm{\tilde{n}}$a and T. Weth, A new look at the fractional Poisson
problem via the logarithmic Laplacian, J. Funct. Anal. 279 (2020), 108732, 50 pp.

\vspace{-0.3cm}

\bibitem{JK81}D. S. Jerison and C. E. Kenig,
The Neumann problem on Lipschitz domains,
Bull. Amer. Math. Soc. (N.S.) 4 (1981), 203--207.

\vspace{-0.3cm}

\bibitem{JK81a}D. S. Jerison and C. E. Kenig, The Dirichlet problem in nonsmooth domains,
Ann. of Math. (2) 113 (1981), 367--382.

\vspace{-0.3cm}

\bibitem{k94}
C. E. Kenig, Harmonic Analysis Techniques for Second Order Elliptic Boundary Value Problems,
CBMS Regional Conf. Ser. in Math. 83,
the American Mathematical Society, Providence, RI, 1994.

\vspace{-0.3cm}

\bibitem{KP93}
C. E. Kenig and J. Pipher,
The Neumann problem for elliptic equations with nonsmooth coefficients,
Invent. Math. 113 (1993), 447--509.

\vspace{-0.3cm}

\bibitem{KP95}
C.~E. Kenig and J. Pipher,
The Neumann problem for elliptic equations with nonsmooth coefficients. II,
Duke Math. J. 81 (1995), 227--250.

\vspace{-0.3cm}

\bibitem{kms15} T. Kuusi, G. Mingione and Y. Sire, Nonlocal self-improving properties,
Anal. PDE 8 (2015), 57-114.

\vspace{-0.3cm}

\bibitem{Le23}
G. Leoni,
A First Course in Fractional Sobolev Spaces,
Graduate Studies in Mathematics, Vol. 229,
American Mathematical Society, Providence, RI, 2023.

\vspace{-0.3cm}

\bibitem{MY26}
W. Ma and S. Yang,
Neumann functions for the fractional Laplace equations on Lipschitz
domains,
Calc. Var. Partial Differential Equations 65 (2026), Paper No. 233, 43 pp.

\vspace{-0.3cm}

\bibitem{McLean00}
W. McLean,
Strongly Elliptic Systems and Boundary Integral Equations,
Cambridge University Press, Cambridge, 2000.

\vspace{-0.3cm}

\bibitem{n21}S. Nowak, Improved Sobolev regularity for linear nonlocal equations with VMO
coefficients, Math. Ann. 385 (2023), 1323--1378.

\vspace{-0.3cm}

\bibitem{n23} S. Nowak, Regularity theory for nonlocal equations with VMO coefficients,
Ann. Inst. H. Poincar\'e C Anal. Non Lin\'eaire 40 (2023), 61--132.

\vspace{-0.3cm}

\bibitem{rs14a}X. Ros-Oton and J. Serra, The Dirichlet problem for the
fractional Laplacian: regularity up to the boundary,
J. Math. Pures Appl. (9) 101 (2014), 275--302.

\vspace{-0.3cm}

\bibitem{RW26}
X. Ros-Oton and M. Weidner,
Regularity for nonlocal equations with local Neumann boundary conditions,
Anal. PDE 19 (2026), 353--411.

\vspace{-0.3cm}

\bibitem{rs14b}X. Ros-Oton and J. Serra, The extremal solution for the fractional Laplacian, Calc.
Var. Partial Differential Equations 50 (2014), 723--750.

\vspace{-0.3cm}

\bibitem{Ry99}
V. S. Rychkov,
On restrictions and extensions of the Besov and Triebel--Lizorkin
spaces with respect to Lipschitz domains,
J. London Math. Soc. (2) 60 (1999), 237--257.

\vspace{-0.3cm}

\bibitem{Sa98}
G. Savar\'e,
Regularity results for elliptic equations in Lipschitz domains,
J. Funct. Anal. 152 (1998), 176--201.

\vspace{-0.3cm}

\bibitem{Shen07}
Z. Shen,
The $L^p$ boundary value problems on Lipschitz domains,
Adv. Math. 216 (2007), 212--254.

\vspace{-0.3cm}

\bibitem{Shen06}
Z. Shen,
Necessary and sufficient conditions for the solvability of the $L^p$
Dirichlet problem on Lipschitz domains,
Math. Ann. 336 (2006), 697--725.

\vspace{-0.3cm}

\bibitem{SY24}
Z. Shi and L. Yao,
New estimates of Rychkov's universal extension operator for Lipschitz
domains and some applications,
Math. Nachr. 297 (2024), 1407--1443.

\vspace{-0.3cm}

\bibitem{s07}L. Silvestre, Regularity of the obstacle problem for a fractional power of the Laplace
operator, Comm. Pure Appl. Math. 60 (2007), 67--112.

\vspace{-0.3cm}

\bibitem{Stein70}
E. M. Stein,
Singular Integrals and Differentiability Properties of Functions,
Princeton Mathematical Series, Vol. 30,
Princeton University Press, Princeton, NJ, 1970.


\vspace{-0.3cm}

\bibitem{Tri02}
H. Triebel,
Function spaces in Lipschitz domains and on Lipschitz manifolds.
Characteristic functions as pointwise multipliers,
Rev. Mat. Complut. 15 (2002), 475--524.

\vspace{-0.3cm}

\bibitem{Ve84}
G. Verchota,
Layer potentials and regularity for the Dirichlet problem for Laplace's
equation in Lipschitz domains,
J. Funct. Anal. 59 (1984), 572--611.

\end{thebibliography}
\end{document}